\documentclass[12pt,oneside,reqno]{amsart}
	\usepackage{txfonts}
	\usepackage{bbm}
	\usepackage{amsmath}
	\usepackage{graphicx}
	\usepackage{mathrsfs}
	\usepackage{stmaryrd}
	\usepackage{enumerate,amsmath,amssymb,amsthm}
	\numberwithin{equation}{section}%公式按章节编号
	\usepackage{amsthm}
	\usepackage{amssymb,amsfonts,mathrsfs}
	\usepackage{enumerate}
	\usepackage{color}
	\usepackage{xcolor}
	\theoremstyle{definition}
	
	\newtheorem{thm}{Theorem}[section]
	\newtheorem{rem}[thm]{Remark}
	\newtheorem{lem}[thm]{Lemma}

	\newtheorem{asn}[thm]{Assumption}

	\theoremstyle{remark}

	\newenvironment{pf}{{\noindent{\sc Proof~} }}{\qed}

	\def\mbf{\mathbf}%加粗实心
	\def\mbb{\mathbb}%空心的字母
	\def\mcl{\mathcal}
	\def\R{{\mathbb R}}
	\def\E{{\mathbb E}}
	\def\S{\mbb S}
	\def\d{\mathrm{d}}
	\def\e{{\mathrm{e}}}

	\def\qed{\hfill$\Box$\medskip}%定理证明结束
	\allowdisplaybreaks

\begin{document}
	
	\title{Stochastic Hamiltonian Type Jump Diffusion Systems with Countable Regimes: Strong Feller Property and Exponential Ergodicity$^{\star}$}
	
	\date{}
	\author{Fubao Xi, Yafei Zhai and Zuozheng Zhang$^{*}$}
	
	\thanks{$^{\star}$  Supported in part by the National Natural Science Foundation of China under Grant No. 12071031.}
	\thanks{$^{*}$ Corresponding author.}
	\thanks{E-mail:  xifb@bit.edu.cn, yafeizhai@bit.edu.cn, zuozhengzhang@mail.bnu.edu.cn}	
	
	\dedicatory{
		School of Mathematics and Statistics,
		Beijing Institute of Technology, \\
		Beijing, 100081,  P.R.China }
	
	\begin{abstract}
		This work focuses on a class of stochastic Hamiltonian type jump diffusion systems with state-dependent switching, in which the switching component has countably infinite many states.
		First, the existence and uniqueness of the underlying processes are
		obtained with the aid of successive construction methods. 
		Then, the Feller property is established by coupling methods.
		Furthermore, the strong Feller property is proved by introducing some auxiliary
		processes and  making use of appropriate Radon-Nikodym derivatives.
		Finally, on the basis of the above results, the exponential ergodicity is
		obtained under the Foster-Lyapunov drift condition. 
	\end{abstract}
	\maketitle
	\rm
	
	%%%%%%%%%%%%%%%%%%%%%%%%%%%%%%%%%%%%%%%%%%%%%%%%%%%%%%%%%%%%%%
	\section{Introduction}
	Consider a Hamiltonian type system subject to random perturbations. More precisely, let $X_1(t)$ and $X_2 (t)$ denote respectively the position and velocity of a particle moving in $\mbb R^d$ at time $t\geq 0$. Suppose that $(X_1(t),X_2(t))$ is governed by the following stochastic differential equation (SDE)
	\begin{equation}\label{eq:baisc1}
		\left\{
		\begin{aligned}
			\text{d}X_{1}(t)&=X_{2}(t)\text{d}t,\\
			\text{d}X_{2}(t)&=b(X_{1}(t),X_{2}(t),\Lambda(t))\text{d}t+\sigma(X_{1}(t),X_{2}(t),\Lambda(t))\text{d}B(t)+\int_{U}c(X_{1}(t-),X_{2}(t-),\Lambda(t),u)N(\text{d}t,\text{d}u),
		\end{aligned}
		\right.
	\end{equation}
	where  $b(x_1,x_2,k)$ and $c(x_1,x_2,k,u)$ are $\mbb R^d$-valued and $\sigma(x_1,x_2,k)$ is $\mbb R^d\times \mbb R^d$-valued for any $x_1,x_2,u\in \mbb R^d$ and $k\in\mbb S:=\{1,2,3,\cdots\}$. Let $(\Omega,\mcl F ,\{\mcl F_t \}_{t\geq 0}, \mbb P)$ be a complete probability space with a filtration $\{\mcl F_t \}_{t\geq 0}$ satisfying the usual conditions (i.e., it is right continuous and $\mcl F_0$ contains all $\mbb P$-null sets), and let $B(t)$ be an $\mcl F_t$-adapted $\mbb R^d$-valued Brownian motion; let $N(\text{d}t,\text{d}u)$
	(corresponding to a random point function $p(t)$) be a Poisson random measure independent of $B(t)$, and let
	$\widetilde{N}(\text{d}t,\text{d}u) = N(\text{d}t,\text{d}u)-\Pi(\text{d}u)\text{d}t$ be the compensated Poisson
	random measure on $[0,\infty)\times U$, where $\Pi(\cdot)$ is a deterministic finite characteristic measure on
	the measurable space
	$(U,\mcl B(U))$. The second component $\Lambda(t)$
	is a right-continuous random jump process with a countably infinite state space $\mbb S$ such that
	\begin{equation}\label{eq:basic2}
		\mathbb{P}\{\Lambda(t+\Delta)=l \mid \Lambda(t)=k,(X_1(t), X_2(t))=(x_1, x_2)\}= \begin{cases}q_{k l}(x_1, x_2) \Delta+o(\Delta), & \text { if }  l \neq k, \\ 1+q_{k k}(x_1, x_2) \Delta+o(\Delta), & \text { if } l=k,\end{cases}
	\end{equation}
	uniformly in $\mbb R^{2d}$, provided $\Delta \downarrow 0$. The matrix $Q(x_1,x_2):=(q_{kl}(x_1,x_2))_{k,l\in\mbb S}$ is the formal generator of the switching process $\Lambda(t)$. 
	
	Note that in addition to the dependence on $x$ and $y$, the functions $b$, $\sigma$ and $c$ also depend on the discrete component $k \in \mathbb{S}$; the motivation for such a formulation will be explained shortly. When they are independent of $k \in \mathbb{S}$, or equivalently in the special case when $\mathbb{S}$ is a singleton set, \eqref{eq:baisc1}  reduces to the usual stochastic Hamiltonian jump diffusion system
	\begin{equation}\label{eq:first}
		\left\{\begin{array}{l}
			\mathrm{d} x_1(t)=x_2(t) \mathrm{d} t, \\
			\mathrm{d} x_2(t)=b(x_1(t), x_2(t)) \mathrm{d} t+\sigma(x_1(t), x_2(t)) \mathrm{d} B(t) + \int_{U}c(x_{1}(t-),x_{2}(t-),u)N(\text{d}t,\text{d}u) .
		\end{array}\right.
	\end{equation}
	By selecting appropriate coefficients for the drift term $b$ and jump function $c$, the model \eqref{eq:first} ($b(x_1,x_2)=-[a(x_1, x_2) x_2+\nabla V(x_1)]$ and $c(x_{1},x_{2},u)=0$) can be interpreted as a Hamiltonian system, with broad applications across various fields of mechanics and physics, including the Duffing, Liénard, and van der Pol equations. In recent decades, there has been growing interest in studying stochastic Hamiltonian systems, with \cite{BAO2022114,DEXHEIMER2022321,Xing2020,2002Stochastic,2001Large, 2010Stochasticzhang} and their references offering insights into \eqref{eq:first} and its variations. From \eqref{eq:baisc1} and \eqref{eq:basic2}, it is evident that $X(t)$ characterizes the jump diffusion behavior, whereas $\Lambda(t)$ describes the switching phenomenon. Additionally, the switching behavior also depends on the jump diffusion component of the state. Thus, we can refer to the process $(X(t), \Lambda(t))$ as a stochastic Hamiltonian type jump diffusion system with state-dependent switching.

	Recently, there has been a significant amount of attention given to the topic of exponential ergodicity for switching diffusion processes and switching jump diffusion processes; see e.g.\cite{1993Meyt3,Shao2015Ergodicity,Shao2016Approximation,2009Asymptotic,2017Exponential,2017feller} and reference therein. Typically, the  diffusion coefficient is non-degenerate. However, in the case of the system described by equations \eqref{eq:baisc1} and \eqref{eq:basic2}, the noise only affects the second component, leading to a degenerate system. Meanwhile, the second component is perturbed by random jumps which describes their state discontinuous changes. Compared with stochastic Hamiltonian systems, the sample paths of stochastic Hamiltonian type jump diffusion systems are not continuous. To the best of our knowledge, there has been limited research on the exponential ergodicity for stochastic Hamiltonian type jump diffusion systems with state-dependent switching. 
	
	%		{\red The system \eqref{eq:baisc1} and \eqref{eq:basic2} belongs to the class of regime-switching diffusions or hybrid diffusions. When $c(x, k, u)=0$, it reduces to a switching diffusion system. When the second component $\Lambda(t)$ is missing, it reduces to the usual Hamiltonian-type system. Such state-dependent switching formulation enables one to describe complex system and their inherent uncertainty and randomness in the environment. However, it adds much difficulty in analysis. Indeed, the dependence between $X(t)$ and $\Lambda(t)$ does complicate problems greatly and it necessitates careful handling of the discrete component $\Lambda(t)$. Owing to their ability to delineate complex system subject to various stochastic perturbations, regime-switching diffusions have received growing attentions recently. Some of the representative works can be found in the two books \cite{2006Stochastic} and \cite{2010Hybrid}. The former one dealt with regime-switching diffusions in which the switching process is a continuous time Markov chain independent of the Brownian motion, whereas the latter one treated processes in which the switching component also depends on the continuous-state component.} 

	In this paper, we first show that under suitable conditions, the solution to the system \eqref{eq:baisc1} and \eqref{eq:basic2} is regular or nonexplosive. To do so, we introduce an another auxiliary Hamiltonian type regime diffusion process without jump $(Y(t),\Lambda^{\prime}(t))$. Through the application of Skorokhod’s representation of jumping processes, we write the discrete event process $\Lambda^{\prime}(t)$ as a stochastic integral with respect to another stationary Poisson point process. Then, we obtain the global existence of a unique strong solution for $(Y(t),\Lambda^{\prime}(t))$ by the interlacing procedure. Finally, by successive construction method (similar to the method of \cite{2009Asymptotic}), we present the existence and uniqueness result for $(X(t),\Lambda(t))$. Moreover, we establish the equivalence between the constructions of the stopping time sequences $\tau^{\prime}_n$ and $\tau^{*}_n$ presented in \cite{2017feller}; see Remark \ref{rem:chong} below for more details.
	
	Then we establish the Feller property for the system \eqref{eq:baisc1} and \eqref{eq:basic2}, making use of the coupling method. The coupling method is widely employed in the study of interacting particle systems and jump diffusion processes, as evidenced by \cite{1989Coupling,1985Interacting,1984Coupling,2006Gradient,Wang2010Regularity} and related references. In this paper, we construct a coupling operator $\mathcal{\widetilde{A}}$ in \eqref{eq:couplingall}. For the coupling process $(\widetilde{X}(t),\widetilde{\Lambda}(t),\widetilde{Z}(t),\widetilde{\Xi}(t))$ 
	generated by $\mathcal{\widetilde{A}}$ starting from $(x, k, z, k)$, it is necessary to handle carefully the first time  when the switching components $\widetilde{\Lambda}$ and $\widetilde{\Xi}$ separate from each other. Additional details can be found in the proof of Theorem \ref{thm:feller}.
	
	To establish the strong Feller property, we follow the proof methodology outlined in \cite{2013Xi,2017feller}. Specifically, we first demonstrate that, under certain conditions, the jump diffusion $X^{(k)}$ in \eqref{steq:yuanshi} possesses the strong Feller property. Next, we establish the strong Feller property for the auxiliary process $(V(t), \Psi(t))$ constructed in equation \eqref{eq:steq:fuzhu}, with the discrete component $\Psi(t)$ being constructed differently from that of (5.2) in \cite{2017feller}. Drawing inspiration from \cite{XI2021856}, we construct a conservation $Q$-matrix $\widehat{Q}=(\widehat{q}_{kl})$ such that 
	\(\widehat{q}_{kl}=\sup_{x\in \R^{2d}}q_{kl}(x) \) for $k\neq l$ and \( \widehat{q}_{kk}=-\sum_{l\neq k} \widehat{q}_{kl}\) for all $k\in \S$, which is more general. Finally, we apply the Radon-Nikodym derivative $M_{T}$ of \eqref{streq:rnde} to establish the strong Feller property for the process $(X(t), \Lambda(t))$.

	Our next focus is on investigating the exponential ergodicity of the system \eqref{eq:baisc1} and \eqref{eq:basic2}. In Section 6, we define the $f$-norm, which is an very strong norm. The well-known total variation norm is actually a special case of the $f$-norm. The $f$-exponential ergodicity was previously studied for jump-diffusion processes with state-dependent switching in \cite{2009Asymptotic} and regime-switching jump diffusion processes with countable regimes in \cite{2017feller}. Assuming that $Q(x, y)$ is irreducible (see Assumption \ref{asn:irreducible} for the precise statement), we further establish the $f$-exponential ergodicity under the Foster-Lyapunov drift condition in Theorem \ref{thm:exerg}.
	
	To facilitate the later presentation, we introduce some frequently used notations here. For $x\in \mbb R^{d}$, $\sigma=\left(\sigma_{i j}\right) \in \mathbb{R}^{d} \times \mathbb{R}^{d}$, define
	$$|x|=\left(\sum_{i=1}^{d}\left|x_{i}\right|^{2}\right)^{1 / 2}, \quad\|\sigma\|=\left(\sum_{i, j=1}^{d}\left|\sigma_{i j}\right|^{2}\right)^{1 / 2}.$$ For $x=(x_1,x_2)$, $y=(y_1,y_2) \in \mathbb{R}^{2 d}$, define a metric $\lambda(\cdot, \cdot)$ on $\mathbb{R}^{2 d} \times \mathbb{S}$ as
	$$
	\lambda((x, m),(y, n)):=|x-y|+\mathbf{1}_{\{m \neq n\}},
	$$
	where $|x|=|(x_1, x_2)|=\sqrt{|x_1|^{2}+|x_2|^{2}}$. Let $\mathcal{B}\left(\mathbb{R}^{2 d} \times \mathbb{S}\right)$ be the Borel $\sigma$-algebra on $\mathbb{R}^{2 d} \times \mathbb{S}$. Then $\left(\mathbb{R}^{2 d} \times \mathbb{S}, \lambda(\cdot, \cdot), \mathcal{B}\left(\mathbb{R}^{2 d} \times \mathbb{S}\right)\right)$ is a locally compact and separable metric space. When there is no special description in what follows, we denote $x=(x_{1},x_{2})$, $y=(y_{1},y_{2})$. In this paper, prior to presenting the proofs of the main theorems, we will provide some essential preparations and warm-up lemmas. The proofs of these lemmas will be deferred to the end of each section.
	
	The paper is organized as follows. In Section 2, we utilize the successive construction method to establish the existence and uniqueness result for stochastic Hamiltonian type jump diffusion systems. Next, in Section 3, we establish the Feller property of stochastic Hamiltonian type jump diffusion systems, making use of the coupling method. In Section 4, we introduce auxiliary processes and apply the Radon–Nikodym derivatives to establish the strong Feller property for stochastic Hamiltonian type jump diffusion systems. Based on these previous results, in Section 5, we further establish the $f$-exponential ergodicity of Hamiltonian type jump diffusion processes under the Foster-Lyapunov drift condition. Finally, in Section 6, we provide a concrete example to illustrate our main results.

	%		Let $C_{c}^{\infty}\left(\mathbb{R}^{2 d} \times \mathbb{S}\right)$ denote the family of functions defined on $\mathbb{R}^{2 d} \times \mathbb{S}$ such that $f(\cdot, k) \in C_{c}^{\infty}\left(\mathbb{R}^{2 d}\right)$ for each $k \in \mathbb{S}$, and $f(x_1, x_2, \cdot)$ is a bounded function on $\mathbb{S}$ for each $(x, y) \in \mathbb{R}^{2 d}$, where $C_{c}^{\infty}\left(\mathbb{R}^{2 d}\right)$ denotes the family of functions defined on $\mathbb{R}^{2 d}$ which are infinitely differentiable and have compact support.
	
	%		As usual, let $C\left([0, \infty), \mathbb{R}^{2 d}\right)$ be the continuous function space endowed with the sup norm topology and $D([0, \infty), \mathbb{S})$ be the càdlàg space endowed with the Skorohod topology. Moreover, let $\Omega:=C\left([0, \infty), \mathbb{R}^{2 d}\right) \times D([0, \infty), \mathbb{S})$ be endowed with the product topology of the sup norm topology on $C\left([0, \infty), \mathbb{R}^{2 d}\right)$ and the Skorohod topology on $D([0, \infty), \mathbb{S})$. Let $\mathcal{F}_{t}$ be the $\sigma$-field generated by the cylindrical sets on $\Omega$ up to time $t$ and set $\mathcal{F}=\bigvee_{t=0}^{\infty} \mathcal{F}_{t}$. Next, 

	\section{Existence and Uniqueness}
	Firstly, we consider the existence and uniqueness of the strong solution to system \eqref{eq:baisc1} and \eqref{eq:basic2}. To do so, we now introduce a family of diffusion processes.
	For each $k\in \mbb S$, let the diffusion process
	$X^{(k),0} (t) = (X^{(k),0}_1 (t),X^{(k),0}_2 (t))$ satisfy the following stochastic differential equation in $\mbb R^{2d}$,
	\begin{equation}\label{steq:yuanshi}
		\left\{\begin{array}{l}
			\text{d} X_{1}^{(k), 0}(t)=X_{2}^{(k),0}(t) \text{d} t, \\
			\text{d} X_{2}^{(k), 0}(t)=b\left(X_{1}^{(k), 0}(t) , X_{2}^{(k), 0}(t), k\right) \text{d} t+\sigma\left(X_{1}^{(k),0}(t), X_{2}^{(k),0}(t), k\right) \text{d} B(t).
		\end{array}\right.
	\end{equation}
	We give a basic assumption for this diffusion process.
	\begin{asn}\label{asn:strongex}
		
		Let $\mcl F_t^{B}$ be a $\sigma$-algebra generated by $\{B(s),s\leq t \}$.	For all $k\in \mbb S$, given a stopping time $\tau$ and an $\mathcal{F}_{\tau}^B$-measurable $\mbb R^{2d}$-valued random variable $x(\tau) $ (depending on $\tau$), there exists a unique strong solution to \eqref{steq:yuanshi} in $[\tau,\infty)$ satisfying $X^{(k),0}(\tau) = x(\tau)$.
	\end{asn}
	\begin{rem}
		Suppose that for all $k\in \mbb S$ and $x$, $y\in \mbb R^{2d}$, there exists a constant $K>0$, such that
		\begin{equation*}
			|b(x,k)-b(y,k)|^2+\|\sigma(x,k)-\sigma(y,k)\|^2\leq K|x-y|^2,\quad|b(x, k)|^{2}+\|\sigma(x, k)\|^{2} \leq K\left(1+|x|^{2}\right),		
		\end{equation*}
		then by Remark 3.10 of \cite{mao}, it is easy to check that Assumption \ref{steq:yuanshi} holds.
	\end{rem}
	Next, we introduce an another auxiliary diffusion process with state-dependent switching. Let 
	$(Y(t),\Lambda^{'}(t))$ satisfy the following stochastic differential equation at $\mbb R^{2d}\times \mathbb{S}$,
	\begin{equation}\label{steq:guocheng}
		\left\{\begin{array}{l}
			\text{d} Y_{1}(t)=Y_{2}(t) \text{d} t, \\
			\text{d} Y_{2}(t)=b\left(Y_{1}(t),Y_{2}(t), \Lambda^{'}(t)\right) \text{d} t+\sigma\left(Y_{1}(t), Y_{2}(t), \Lambda^{'}(t) \right) \text{d} B(t),
		\end{array}\right.
	\end{equation}
	\begin{equation}\label{steq:switch}
		\mathbb{P}\{\Lambda^{'}(t+\Delta)=l \mid \Lambda^{'}(t)=k,(Y_1(t), Y_2(t))=(y_1, y_2)\}= \begin{cases}q_{k l}(y_1, y_2) \Delta+o(\Delta), & \text { if } l \neq k, \\ 1+q_{k k}(y_1, y_2) \Delta+o(\Delta), & \text { if } l=k\end{cases}
	\end{equation}
	uniformly in $\mbb R^{2d}$, provided $\triangle \downarrow 0$. 
	\begin{asn}\label{asn:zysl}
		Assume that for all $(x,k)\in \mbb R^{2d}\times \mbb S$, we have 
		\begin{equation}\label{eq:yizhi}
			q_k(x):=-q_{kk}(x)=\sum_{l\in \mbb S\backslash\{k\}}q_{kl}(x)\leq H
		\end{equation}
		for some constant $H>0$.
	\end{asn}
	\begin{lem}\label{prob:stong}
		Under Assumptions \ref{asn:strongex} and \ref{asn:zysl}, for each initial data $(y, k)$, system \eqref{steq:guocheng} and \eqref{steq:switch} has a unique strong solution $(Y(t), \Lambda^{'}(t))$ with $(Y(0), \Lambda^{'}(0))=(y, k)$, and $(Y(t), \Lambda^{'}(t))$ is non-explosive.
	\end{lem}
	\begin{thm}\label{prob:sstong}
		Under Assumptions \ref{asn:strongex} and \ref{asn:zysl}, for each initial data $(x, k)$, system \eqref{eq:baisc1} and \eqref{eq:basic2} has a unique strong solution $(X(t), \Lambda(t))$ with $(X(0), \Lambda(0))=(x, k)$, and $(X(t), \Lambda(t))$ is non-explosive.
	\end{thm}
	
	\noindent\textbf{\textit{Proof.}} Let $\sigma_{1}<\sigma_{2}<\cdots<\sigma_{n}<\cdots$ be the enumeration of all elements in the domain $D_{p}$ of the point process $p(t)$ corresponding to the above Poisson random measure $N(\d t, \d u)$. It is easy to see that $\sigma_{n}$ is an $\mathcal{F}_{t}$-stopping time for each $n$ and $\left\{\sigma_{n+1}-\sigma_{n}\right\}_{n \geq 0}$ is a sequence of independent and identically distributed exponential random variables, where $\sigma_{0}=0$. Moreover, we have that $\lim_{n \rightarrow \infty} \sigma_{n}=+\infty$ almost surely since the characteristic measure $\Pi(\cdot)$ is a finite measure on the measurable space $(U,\mcl B(U))$. Let us denote the successive switching stopping times of the second component $\Lambda(t)$ by
	$$
	\tau_0 \equiv 0, \quad \tau_n=\inf \left\{t: t>\tau_{n-1}, \Lambda(t) \neq \Lambda\left(\tau_{n-1}\right)\right\}, \quad n \geq 1.
	$$
	Here, we will establish the existence and uniqueness of a solution $(X(t), \Lambda(t))$ to system \eqref{eq:baisc1} and \eqref{eq:basic2}. Firstly, let us consider it in the time interval $\left[0, \sigma_1\right]$. For any $t \in\left[0, \sigma_1\right)$ and any path $\{(X(s), \Lambda(s)): 0 \leq s \leq t\}$, we always have
	$$
	\int_0^t \int_{U} c(X(s-), \Lambda(s-), u) N(\mathrm{d}s,\mathrm{d}u) \equiv 0 .
	$$
	Hence, on the interval $\left[0, \sigma_1\right)$, we can consider the system \eqref{steq:guocheng} and \eqref{steq:switch} instead of system \eqref{eq:baisc1} and \eqref{eq:basic2}. By Lemma \ref{prob:stong}, there exists a unique strong solution $\left(Y(t), \Lambda^{\prime}(t)\right)$ to the system \eqref{steq:guocheng} and \eqref{steq:switch} such that $\left(Y(0), \Lambda^{\prime}(0)\right)=(x, k)$. Since $\sigma_1$ does not coincide with any of $\left\{\tau_n: n \geq 1\right\}$ almost surely by Proposition 2.3 of \cite{2009Asymptotic}, we can construct the strong Markov process $(X(t), \Lambda(t))$ as follows. For any given initial condition $(X(0), \Lambda(0))=(x, k)$, on the time interval $\left[0, \sigma_1\right]$, set
	$$
	(X_1(t),X_2(t))=\left\{\begin{array}{ll}
		(Y_1(t),Y_2(t)), & 0 \leq t<\sigma_1, \\
		\Big(Y_1(\sigma_1-  ),Y_2(\sigma_1-)+c\big(Y_1(\sigma_1-),Y_2\left(\sigma_1-\right), \Lambda^{\prime}\left(\sigma_1-\right), p\left(\sigma_1\right)\big)\Big), & t=\sigma_1,
	\end{array}\right.
	$$
	and
	$$
	\Lambda(t)=\Lambda^{\prime}(t), \quad 0 \leq t \leq \sigma_1,
	$$
	where $p(t)$ is the random point process corresponding to $N(\mathrm{d} t, \mathrm{d} u)$. Next, as previously done in Chapter IV, Section 9 of \cite{1983Stochastic},  $\widetilde{x}=X\left(\sigma_1\right), \widetilde{k}=\Lambda\left(\sigma_1\right), \widetilde{B}(t)=B\left(t+\sigma_1\right)-B\left(\sigma_1\right)$ and $\widetilde{p}(t)=p\left(t+\sigma_1\right)$. Recall that $\sigma_2$ does not coincide with any of $\left\{\tau_n: n \geq 1\right\}$ almost surely. Similarly, we can determine the process $(\widetilde{X}(t), \widetilde{\Lambda}(t))$ on the time interval $\left[0, \sigma_2-\sigma_1\right]$ with respect to $(\widetilde{x}, \widetilde{k})$ as above. Then, define
	$$
	(X(t), \Lambda(t))=\left(\widetilde{X}\left(t-\sigma_1\right), \widetilde{\Lambda}\left(t-\sigma_1\right)\right), \quad t \in\left[\sigma_1, \sigma_2\right] .
	$$
	Continuing this procedure successively, $(X(t), \Lambda(t))$ is determined uniquely on the time interval $\left[0, \sigma_n\right]$ for every $n$. Thus, $(X(t), \Lambda(t))$ is determined globally due to $\lim_{n \rightarrow \infty} \sigma_n=+\infty$. Therefore, we have proved the existence and uniqueness of the strong solution to system \eqref{eq:baisc1} and \eqref{eq:basic2}.\hfill $\square$\par
	
	In what follows, we shall give the proof of Lemma \ref{prob:stong}. Note that the evolution of the discrete component $\Lambda^{'}(t)$ can be represented as a stochastic integral with respect to a stationary Poisson point process through the application of Skorokhod’s representation. Indeed, for each $k \in \mbb S$ and $y \in \mathbb{R}^{2 d}$ , let
	\[
	\Delta_{k 1}(y)=\left[0, q_{k 1}(y)\right), \Delta_{k l}(y)=\left[\sum_{j=1,j\neq k}^{l-1}q_{kj}(y),\sum_{j=1,j\neq k}^{l}q_{kj}(y)\right),l>1,l\neq k.
	\]
	Note that for each $k \in \mbb S$ and $y \in \mathbb{R}^{2 d}$, $\{\Delta_{kl}(y):l\in \mbb S \}$ are disjoint intervals, and the length of the interval $\Delta_{kl}(y)$ is equal to $q_{kl}(y)$, which is bounded above by $H$ thanks to Assumption \ref{asn:zysl}. We then define a function $h: \mathbb{R}^{2d} \times \mathbb{S} \times\left[0, H\right]$ by
	\[
	h(y, k, u)=\sum_{l \in \mbb S}(l-k) \mathbf{1}_{\triangle_{k l}(y)}(u).
	\]
	That is, for each $k \in \mathbb{S}$, if $u \in \triangle_{k l}(y)$, $h(y, k, u)=l-k$; otherwise $h(y, k, u)=0$. Hence, \eqref{steq:switch} is equivalent to
	\begin{equation}\label{eq:skorohod}
		\text{d}\Lambda^{'}(t)=\int_{\left[0, H\right]} h(Y(t), \Lambda^{'}(t-), u) N_1(\text{d} t, \text{d} u),
	\end{equation}
	where $N_1(\d t, \d u)$ is a Poisson random measure (corresponding to a random point process $p_1(t)$) with the Lebesgue measure on $\left[0, H\right]$ as its characteristic measure. As mentioned in \cite{1999Stability}, the process $(Y(t), \Lambda^{'}(t))$ can be thought of as a solution to system \eqref{steq:guocheng} and \eqref{eq:skorohod} with the driving forces being the Brownian motion $B(\cdot)$ and the Poisson point process $N_1(\cdot, \cdot)$ which is independent of $B(\cdot)$. 
	%		Denote the natural filtration by $\mathcal{F}^{\prime}_{t}:=\sigma\{(Y(s), \Lambda^{'}(s)): s \leq t\}$. Without loss of generality, assume that the filtration $\left\{\mathcal{F}^{\prime}_{t}\right\}_{t \geq 0}$ satisfies the usual condition (i.e., it is right continuous and $\mathcal{F}^{\prime}_{0}$ contains all $P$-null sets).
	
	\noindent\textbf{\textit{Proof of Lemma \ref{prob:stong}.}} Note that Assumption \ref{asn:strongex} guarantees the existence and uniqueness of strong solutions to the following diffusion: for each $i \in \mathbb{S}$,
	\begin{equation}\label{steq:baiscc}
		\left\{\begin{array}{l}
			\text{d} Y_{1}(t)=Y_{2}(t)\text{d} t, \\
			\text{d} Y_{2}(t)=b\left(Y_{1}(t),Y_{2}(t), i\right) \text{d} t+\sigma\left(Y_{1}(t), Y_{2}(t), i \right) \text{d} B(t).
		\end{array}\right.
	\end{equation}
	Given a stopping time $\tau$ and an $\mathcal{F}_\tau^B$-measurable $\mathbb{R}^{2d}$-valued random variable $y=y(\tau)$, there exists a unique strong solution to \eqref{steq:baiscc} in $[\tau, \infty)$ satisfying $Y(\tau)=y(\tau)$ . We can now construct the solution to system \eqref{steq:guocheng} and \eqref{eq:skorohod} with initial data $\left(y, k\right)$ by the interlacing procedure similar to Chapter 5 of \cite{applebaum_2009}. Let $k_0=k,\tau^{\prime}_0=0$ for convenience, and let $Y^{(0)}(t)$, $t \geq 0$ be the solution with initial data $y=(y_1,y_2)$ to
	\begin{equation}\label{eq:chifan}
		\left\{\begin{array}{l}
			\text{d} Y^{(0)}_{1}(t)=Y^{(0)}_{2}(t) \text{d} t, \\
			\text{d} Y^{(0)}_{2}(t)=b\left(Y^{(0)}_{1}(t),Y^{(0)}_{2}(t), k_0\right) \text{d} t+\sigma\left(Y^{(0)}_{1}(t), Y^{(0)}_{2}(t), k_0 \right) \text{d} B(t).
		\end{array}\right.
	\end{equation}
	Let
	\[
	\tau^{\prime}_1:=\inf \left\{t>0: \int_{0}^t \int_{\left[0, H\right]} h(Y^{(0)}(s),k_0, u) N_1(\text{d} s, \text{d} u) \neq 0\right\} \text { and } 
	\]
	\[
	k_1:=k_0+\int_{0}^{\tau^{\prime}_1} \int_{\left[0, H\right]} h(Y^{(0)}(s),k_0, u) N_1(\text{d} s, \text{d} u),
	\]
	and let $Y^{(1)}(t)$, $t \geq \tau^{\prime}_1$ be the solution with $Y^{(1)}(\tau^{\prime}_1)=Y^{(0)}(\tau^{\prime}_1)$ to 
	\begin{equation}
		\left\{\begin{array}{l}
			\text{d} Y^{(1)}_{1}(t)=Y^{(1)}_{2}(t) \text{d} t ,\\
			\text{d} Y^{(1)}_{2}(t)=b\left(Y^{(1)}_{1}(t),Y^{(1)}_{2}(t), k_1\right) \text{d} t+\sigma\left(Y^{(1)}_{1}(t), Y^{(1)}_{2}(t), k_1 \right) \text{d} B(t).
		\end{array}\right.
	\end{equation}
	Let
	\[
	\tau^{\prime}_2:=\inf \left\{t>\tau^{\prime}_1: \int_{\tau^{\prime}_1}^t \int_{\left[0, H\right]} h(Y^{(1)}(s),k_1, u) N_1(\text{d} s, \text{d} u) \neq 0\right\} \text { and } 
	\]
	\[
	k_2:=k_1+\int_{\tau^{\prime}_1}^{\tau^{\prime}_2} \int_{\left[0, H\right]} h(Y^{(1)}(s),k_1, u) N_1(\text{d} s, \text{d} u).
	\]
	Continuing this procedure, let $\tau^{\prime}_{\infty}=\lim _{k \rightarrow \infty} \tau^{\prime}_k$ and set
	$$
	Y(t)=Y^{(i)}(t), \quad\Lambda^{\prime}(t)=k_i \quad\text { if } \tau^{\prime}_i \leq t<\tau^{\prime}_{i+1} .
	$$
	Clearly, $Y(t)$ satisfies that for every $t \geq 0$,
	\begin{equation*}
		\left\{\begin{array}{l}
			Y_{1}(t\wedge \tau^{\prime}_i)=y_1+\int_0^{t \wedge \tau^{\prime}_i}Y_{2}(s) \text{d} s, \\
			Y_{2}(t\wedge \tau^{\prime}_i)=y_2+\int_0^{t \wedge \tau^{\prime}_i}b\left(Y_{1}(s),Y_{2}(s), \Lambda^{'}(s)\right) \text{d} s+\int_0^{t \wedge \tau^{\prime}_i}\sigma\left(Y_{1}(s), Y_{2}(s), \Lambda^{'}(s) \right) \text{d} B(s),
		\end{array}\right.
	\end{equation*}
	\begin{equation*}
		\Lambda^{'}(t\wedge\tau^{\prime}_i)=k_0+\int_0^{t \wedge \tau^{\prime}_i}\int_{\left[0, H\right]} h(Y(s),\Lambda^{'}(s-), u) N_1(\text{d} s, \text{d} u),
	\end{equation*}
	To show that $(Y(t),\Lambda^{'}(t))$ is a global solution, we only need  to prove that $\tau^{\prime}_{\infty}=\infty$ a.s. By Assumption \ref{asn:zysl}, for any $T>0$,
	$$
	\begin{aligned}
		\mathbb{P}\left\{\tau^{\prime}_i \leq T\right\} &=\mathbb{P}\left\{\int_0^{T \wedge \tau^{\prime}_i} \int_{H} \mathbf{1}_{\left\{u \in\left[0, q_{\Lambda^{\prime}(s-)}\left(Y(s)\right)\right)\right\}} N_1(\d s, \d u)=i\right\} \\
		& \leq \mathbb{P}\left\{\int_0^{T \wedge \tau^{\prime}_i} \int_{H} \mathbf{1}_{\{u \in[0, H)\}} N_1(\d s, \d u) \geq i\right\} \\
		& \leq \mathbb{P}\left\{\int_0^T \int_{H} \mathbf{1}_{\{u \in[0, H)\}} N_1(\d s, \d u) \geq i\right\} \\
		&=\sum_{l=i}^{\infty} \e^{-H T} \frac{(H T)^l}{l !} .
	\end{aligned}
	$$
	It follows that $\mathbb{P}\left\{\tau^{\prime}_i \leq T\right\} \rightarrow 0$ as $i \rightarrow \infty$. As a result,  $\tau^{\prime}_{\infty}=\infty$ a.s. By this construction, it can be seen that $Y(t)$ is continuous and $\Lambda^{\prime}(t)$ is cadlag a.s. The uniqueness of $(Y(t), \Lambda^{\prime}(t))$ follows from the uniqueness of $Y^{(i)}(t)$ on $\left[\tau^{\prime}_i, \tau^{\prime}_{i+1}\right]$ and the uniqueness of $k_i$ defined by
	$$
	k_{i}=k_{i-1}+\int_{\tau^{\prime}_{i-1}}^{\tau^{\prime}_{i}} \int_{\left[0, H\right]} h(Y^{(i-1)}(s),k_{i-1}, u) N_1(\text{d} s, \text{d} u).
	$$
	The proof is complete.\hfill $\square$\par		
	\begin{rem}\label{rem:chong}
		There is an equivalent way to construct stopping time sequence  in \cite{2017feller}. Let us provide a brief description of this construction method as follows. Let $\{\xi_n\}$ be a sequence of independent mean $1$ exponential random variables on $(\Omega,\mathcal{F},\{\mathcal{F}_t\}_{t\geq 0},\mbb P)$ independent of $B(t)$. Let 
		\[
		\tau^{*}_1=\theta_1:=\inf\left\{t\geq 0:\int_{0}^t q_k(Y^{(0)}(s))\d s>\xi_1\right\}.
		\]
		Then,
		\[
		\mbb P\{\tau_1^*>t|\mathcal{F}_t^B\}=\mbb P\left\{\xi_1\geq \int_{0}^{t}q_k(Y^{(0)}(s))\d s\big|\mathcal{F}_t^B\right\}=\exp\left\{-\int_{0}^{t}q_k(Y^{(0)}(s))\d s\right\}.
		\]
		Moreover, we define \( \Lambda^{\prime}\left(\tau_{1}^{*}\right)  \) according to the probability distribution:
		\[
		\mathbb{P}\left\{\Lambda^{\prime}\left(\tau_{1}^{*}\right)=l \big|\mathcal{F}_{\tau_{1}^{*}-}^B\right\}=\frac{q_{k l}\left(Y^{(0)}\left(\tau_{1}^{*}-\right)\right)}{q_{k}\left(Y^{(0)}\left(\tau_{1}^{*}-\right)\right)}\left(1-\delta_{k l}\right) \mathbf{1}_{\left\{q_{k}\left(Y^{(0)}\left(\tau_{1}^{*}-\right)\right)>0\right\}}+\delta_{k l} \mathbf{1}_{\left\{q_{k}\left(Y^{(0)}\left(\tau_{1}^{*}-\right)\right)=0\right\}} ,
		\]
		where $Y^{(0)}(t)$, $t \geq 0$ is defined in \eqref{eq:chifan}. Let $Y^{(1)}(t)$, $t \geq 0$ be the solution with $Y^{(1)}(0)=Y^{(0)}(\tau^{*}_1)$ to 
		\begin{equation}
			\left\{\begin{array}{l}
				\text{d} Y^{(1)}_{1}(t)=Y^{(1)}_{2}(t) \text{d} t, \\
				\text{d} Y^{(1)}_{2}(t)=b\left(Y^{(1)}_{1}(t),Y^{(1)}_{2}(t), \Lambda^{\prime}\left(\tau_{1}^{*}\right)\right) \text{d} t+\sigma\left(Y^{(1)}_{1}(t), Y^{(1)}_{2}(t), \Lambda^{\prime}\left(\tau_{1}^{*}\right) \right) \text{d} B(t).
			\end{array}\right.
		\end{equation}
		Let
		\[
		\theta_2:=\inf \left\{t \geq 0: \int_{0}^{t} q_{\Lambda^{\prime}\left(\tau_1^*\right)}\left(Y^{\left(1\right)}(s)\right) \mathrm{d} s>\xi_{2}\right\},\qquad \tau_2^*:=\tau_1^*+\theta_2.
		\]
		Then,
		\[
		\begin{aligned}
			\mathbb{P}\left\{\theta_2>t \big| \mathcal{F}_{\tau_1^*+t}^B\right\} & =\mathbb{P}\left\{\xi_{2} \geq \int_{0}^{t} q_{\Lambda^{\prime}\left(\tau_1^*\right)}\left(Y^{\left(1\right)}(s)\right) \mathrm{d} s \big| \mathcal{F}_{\tau_1^*+t}^B\right\} 
			& =\exp \left\{-\int_{0}^{t} q_{\Lambda^{\prime}\left(\tau_1^*\right)}\left(Y^{\left(1\right)}(s)\right) \mathrm{d} s\right\} .
		\end{aligned}
		\]
		Similarly, we define \( \Lambda^{\prime}\left(\tau_{2}^{*}\right) \):
		\[
		\begin{aligned}
			&\mathbb{P}  \left\{\Lambda^{\prime}\left(\tau_{2}^*\right)=l \mid \mathcal{F}_{\tau_{2}^*-}\right\} \\
			 =&\frac{q_{\Lambda^{\prime}(\tau_{1}^*), l}\left(Y^{(1)}\left(\tau_{2}^*-\right)\right)}{q_{\Lambda^{\prime}\left(\tau_{1}^*\right)}\left(Y^{(1)}\left(\tau_{2}^*-\right)\right)}\left(1-\delta_{\Lambda^{\prime}\left(\tau_{1}^*\right), l}\right) \mathbf{1}_{\left\{q_{\Lambda^{\prime}\left(\tau_{1}^*\right)}\left(Y^{(1)}\left(\tau_{2}^*-\right)\right)>0\right\}}+\delta_{\Lambda^{\prime}\left(\tau_{1}^*\right), l} \mathbf{1}_{\left\{q_{\Lambda^{\prime}\left(\tau_{1}^*\right)}\left(Y^{(1)}\left(\tau_{2}^*-\right)\right)=0\right\}} .
		\end{aligned}
		\]
		Continuing this procedure, let $\tau^{*}_{\infty}=\lim _{n \rightarrow \infty} \tau^{*}_n$ and set
		$$
		Y(t)=Y^{(i)}(t-\tau^{*}_n), \quad\Lambda^{\prime}(t)=\Lambda^{\prime}(\tau^{*}_n ) \quad\text { if } \tau^{*}_n \leq t<\tau^{*}_{n+1} ,
		$$
		and
		\[
		\begin{aligned}
			&\mathbb{P}  \left\{\Lambda^{\prime}\left(\tau_{n+1}^*\right)=l \mid \mathcal{F}_{\tau_{n+1}^*-}\right\} \\
			=& \frac{q_{\Lambda^{\prime}(\tau_{n}^*), l}\left(Y^{(n)}\left(\tau_{n+1}^*-\right)\right)}{q_{\Lambda^{\prime}\left(\tau_{n}^*\right)}\left(Y^{(n)}\left(\tau_{n+1}^*-\right)\right)}\left(1-\delta_{\Lambda^{\prime}\left(\tau_{n}^*\right), l}\right) \mathbf{1}_{\left\{q_{\Lambda^{\prime}\left(\tau_{n}^*\right)}\left(Y^{(n)}\left(\tau_{n+1}^*-\right)\right)>0\right\}}+\delta_{\Lambda^{\prime}\left(\tau_{n}^*\right), l} \mathbf{1}_{\left\{q_{\Lambda^{\prime}\left(\tau_{n}^*\right)}\left(Y^{(n)}\left(\tau_{n+1}^*-\right)\right)=0\right\}} .
		\end{aligned}
		\]
		By Assumption \ref{asn:zysl}, we have \( \mathbb{P}\left\{\theta_{n}>t\right\} \geq \e^{-H t} \) for all \( n \in \mathbb{N} \) and \( t>0 \). Hence
		\[
		\begin{aligned}
		&	\mathbb{P}\left\{\tau^*_{\infty}=\infty\right\}  \geq \mathbb{P}\left\{\left\{\theta_{k}>t\right\} \text { i.o. }\right\}=\mathbb{P}\left\{\bigcap_{m=1}^{\infty} \bigcup_{k=m}^{\infty}\left\{\theta_{k}>t\right\}\right\} \\
			=& \lim _{m \rightarrow \infty} \mathbb{P}\left\{\bigcup_{k=m}^{\infty}\left\{\theta_{k}>t\right\}\right\} \geq \limsup _{m \rightarrow \infty} \mathbb{P}\left\{\theta_{m}>t\right\} \geq \e^{-H t} .
		\end{aligned}
		\]
		Letting \( t \downarrow 0 \) yields that \( \mathbb{P}\left\{\tau^*_{\infty}=\infty\right\}=1 \). This shows that $(Y(t),\Lambda^{'}(t))$ is a global solution.

		Next, we shall show that $\tau^{\prime}_n$ and $\tau^{*}_n$ have the same conditional distribution. Due to the construction of $\tau^{\prime}_n$ and $\tau^{*}_n$, it suffices to show that $\tau^{\prime}_1$ and $\tau^{*}_1$ have the same conditional distribution. Indeed, let $V(j)=1$ if $j=k$ and $V(j)=0$ if $j\neq k$. Applying the generalized Itô formula to $V$, we have
		\[
		\begin{aligned}
			&\mbb P\{\tau_1^{\prime}>t|\mathcal{F}_T^B\}=\mbb E\left[\mbf1_{\{\tau_1^{\prime}>t\}}|\mathcal{F}_T^B\right]=\mbb E\left[V(\Lambda^{\prime}(\tau_1^{\prime}\wedge t))|\mathcal{F}_T^B \right]\\
			=&1+\mbb E[\int_{0}^{\tau_1^{\prime}\wedge t}\sum_{j\in \S}q_{kj}(Y^{(0)}(s))V(j)\d s|\mathcal{F}_T^B ]=1-\E [\int_{0}^{\tau_1^{\prime}\wedge t}q_k(Y^{(0)}(s))\d s|\mathcal{F}_T^B ]\\
			=&1-\E[\int_{0}^t q_k(Y^{(0)}(s))\mbf1_{\{ \tau_1^{\prime}>s\}}\d s |\mathcal{F}_T^B]=1-\int_{0}^t q_k(Y^{(0)}(s))\E[\mbf1_{\{ \tau_1^{\prime}>s\}}|\mathcal{F}_T^B]\d s.
		\end{aligned}
		\]
		Hence 
		\[
		\frac{\d}{\d t}\mbb E\left[\mbf1_{\{\tau_1^{\prime}>t\}}|\mathcal{F}_T^B\right]=-q_k(Y^{(0)}(t))\E[\mbf1_{\{ \tau_1^{\prime}>t\}}|\mathcal{F}_T^B],
		\]
		which, together with $\mbb E\left[1_{\{\tau_1^{\prime}>0\}}|\mathcal{F}_T^B\right]=1$, implys for any $0\leq t\leq T$,
		\[
		\mbb P\{\tau_1^{\prime}>t|\mathcal{F}_T^B\}=\exp\left\{-\int_{0}^{t}q_k(Y^{(0)}(s))\d s\right\}.
		\]
		Letting $t=T$ in the above equation implies the desired assertion.
	\end{rem}
	\section{Feller Property}
	In the previous section, we show that there exists a unique strong solution for system \eqref{eq:baisc1} and \eqref{eq:basic2}. In this section, we shall show that $(X(t),\Lambda(t))$ admits Feller property. Before we proceed, let us give some additional assumptions.
	\begin{asn}\label{asn:basic}
		Suppose that for any $k\in \mbb S$, there exists a positive constant $L_k$ such that, for all $x$, $y\in \mbb R^{2d}$,
		\begin{equation}\label{eq:asnbasic1}
			|b(x, k)|^{2}+\|\sigma(x, k)\|^{2}+\int_{U}|c(x, k, u)|^{2} \Pi(\mathrm{d} u) \leq L_k\left(1+|x|^{2}\right)
		\end{equation}
		\begin{equation}\label{eq:asnbasic2} 
			|b(x,k)-b(y,k)|^2+\|\sigma(x,k)-\sigma(y,k)\|^2+\int_{U}|c(x,k,u)-c(y,k,u)|^2\Pi(\d u)\leq L_k|x-y|^2.			
		\end{equation}
	\end{asn}
	\begin{asn}\label{asn:feller}
		Suppose that for any $k\in \mbb S$, there exists a positive constant $M_k$ such that, for all $x$, $y\in \mbb R^{2d}$,
		\begin{equation}\label{eq:asnfeller2}
			\sum_{l \in \mathbb{S} \backslash\{k\}}\left|q_{k l}(x)-q_{k l}(y)\right| \leq M_k|x-y|.
		\end{equation}
	\end{asn}
	\begin{rem}
		It is clear that Assumption \ref{asn:strongex} holds under Assumption \ref{asn:basic}. Moreover, note that the Lipschitz cosntants of $b(\cdot,k)$, $\sigma(\cdot,k)$ and $c(\cdot,k,u)$ all depend on $k$. If Lipschitz constants are independent of $k$ (see Assumption \ref{asn:basic2}), the uniform boundedness for $q_k(x)$ (see \eqref{eq:yizhi}) will be no longer required. In this case, Assumption \ref{asn:zysl} can be replaced by the following assumption:
		
		\noindent\textbf{Assumpiton 2.3$^{\prime}$}	Assume that for any $k\in \mbb S$ and $M>0$, there exists some constant $H>0$ such that
		\begin{equation*}
			\sup_{|x|\leq M} q_k(x)=\sup_{|x|\leq M}\sum_{l\in \mbb S\backslash\{k\}}q_{kl}(x)\leq H.
		\end{equation*}
		Further, with a uniform Lipschitz constant and the above assumption, the existence and uniqueness of solutions to system \eqref{eq:baisc1} and \eqref{eq:basic2} can be obtained; see Theorem 3.3 in \cite{NG2016} for more details.
	\end{rem}
	It is well-known that $(X(t), \Lambda(t))$ can also be associated with an appropriate generator. Let $\langle\cdot, \cdot\rangle$ and $\nabla$ denote the inner product and the gradient operator in $\mathbb{R}^{2d}$, respectively. As usual, let $C_{c}^{2}\left(\mathbb{R}^{2d}\right)$ denote the family of all functions on $\mathbb{R}^{2d}$ which are twice continuously differentiable and have compact support. For each $k \in \mbb S$, and for any function $f(\cdot, k) \in C_{c}^{2}\left(\mathbb{R}^{2d}\right)$, define an operator $\mathcal{A}$ as follows:
	\begin{equation}\label{eq:alloperator}
		\mathcal{A}f(x,k):=\mathcal{L}_{k}f(x,k)+\Omega(k) f(x,k)+Q(x)f(x,k).
	\end{equation}
	Here operators $\mathcal{L}_k$, $\Omega(k)$ and $Q(x)$ are further defined as follows, for all $x\in \mbb R^{2d}$ and $k\in \mbb S$,
	\begin{equation}
		\mathcal{L}_{k}f(x,k)=\langle\nabla_{x_{1}}f(x,k),x_{2}\rangle+\langle b(x,k),\nabla_{x_{2}}f(x,k)\rangle+\frac{1}{2}{\rm tr}(a(x,k)\nabla_{x_{2}}^2f(x,k)),
	\end{equation}
	\begin{equation}\label{eq:operatorOMega}
		\Omega(k)f(x,k)=\int_{U}[f(x_{1},x_{2}+c(x,k,u),k)-f(x_1,x_2,k)]\Pi(\text{d}u),
	\end{equation}
	\begin{equation}\label{eq:operatorQ}
		Q(x)f(x,k)=\sum_{l\in \S}q_{kl}(x)(f(x,l)-f(x,k)),
	\end{equation}
	where $a(x,k)=\sigma(x,k)\sigma(x,k)^\prime$.
	
	In this section, we will further establish the Feller property of $(X(t),\Lambda(t))$ using coupling methods. To this end, we will first construct a coupling operator $\widetilde{\mathcal{A}}$ for $\mathcal{A}$. For $x$, $y\in \mbb R^{2d}$ and $k$, $l\in \mbb S$, we get
	\[
	a(x,k,y,l)=\begin{pmatrix}
		a(x,k)&\sigma(x,k)\sigma(y,l)^\prime\\
		\sigma(y,l)\sigma(x,k)^\prime&a(y,l)
	\end{pmatrix},
	\qquad
	b(x,k,y,l)=\begin{pmatrix}
		b(x,k)\\
		b(y,l)
	\end{pmatrix}.
	\]
	Next, for $f(x,k,y,l)\in C_{c}^{2}(\mbb R^{2d}\times \mbb S\times \mbb R^{2d}\times \mbb S)$, we define
	\begin{equation}
		\begin{aligned}
			\begin{aligned}
				\widetilde{\Omega}_{d}f(x,k,y,l)=&\langle\nabla_{x_{1}}f(x,k,y,l),x_{2}\rangle+\langle\nabla_{y_{1}}f(x,k,y,l),y_{2}\rangle\\
				&+\langle b(x,k,y,l),Df(x,k,y,l)\rangle+\frac{1}{2}{\rm tr}(a(x,k,y,l)D^2f(x,k,y,l)),
			\end{aligned}
		\end{aligned}
	\end{equation}
	where in the above, $Df(x,k,y,l)$ represents the gradient of $f$ with respect to the variables $x_{2}$ and $y_{2}$, that is, $Df(x,k,y,l)=(\nabla_{x_{2}}f(x,k,y,l),\nabla_{y_{2}}f(x,k,y,l))^\prime$. Likewise, $D^{2}f(x,k,y,l)$ denotes the Hessian of $f$ with respect to the variables $x_{2}$ and $y_{2}$. Let us also define for $f(x,k,y,l)\in C_{c}^{2}(\mbb R^{2d}\times \mbb S\times \mbb R^{2d}\times \mbb S)$,
	%\begin{small}
	\begin{equation}
		\begin{aligned}
			\widetilde{\Omega}_{j}f(x,k,y,l)=&\int_{U}[f(x_1,x_2+c(x,k,u),k,y_1,y_2+c(y,l,u),l)-f(x_1,x_2,k,y_1,y_2,l)]\Pi(\d u),
		\end{aligned}
	\end{equation}
	%\end{small}
	which is a coupling of the jump part in the generator $\Omega(k)$ defined in \eqref{eq:operatorOMega}. 
	
	Then we define the basic coupling \cite{2004From}  for the $q$-matrices $Q(x)$ and $Q(y)$. For any $f(x,k,y,l)\in C_{c}^{2}(\mbb R^{2d}\times \mbb S\times \mbb R^{2d}\times \mbb S)$, we define
	\begin{equation}\label{eq:couoperQ}
		\begin{aligned}
			&\widetilde{\Omega}_{s}f(x,k,y,l)\\
			&=\sum_{i\in \S}[q_{ki}(x)-q_{li}(y)]^+(f(x,i,y,l)-f(x,k,y,l))+\sum_{i\in \S}[q_{li}(y)-q_{ki}(x)]^+(f(x,k,y,i)-f(x,k,y,l))\\
			&\quad+\sum_{i\in \S}[q_{ki}(x)\wedge q_{li}(y)](f(x,i,y,i)-f(x,k,y,l)).
		\end{aligned}
	\end{equation}
	It is easy to verify that $\widetilde{Q}(x,y)$ defined in \eqref{eq:couoperQ} is a coupling for $Q(x)$ defined in \eqref{eq:operatorQ}.
	Finally, the coupling operator of $\mathcal{A}$ in \eqref{eq:alloperator} can be written as
	\begin{equation}\label{eq:couplingall}
		\widetilde{\mathcal{A}}f(x,k,y,l):=[\widetilde{\Omega}_{d}+\widetilde{\Omega}_{j}+\widetilde{\Omega}_{s}]f(x,k,y,l).
	\end{equation}
	In fact, we can verify directly that for any $f(x,k,y,l)=g(x,k)\in C_{c}^{2}(\mbb R^{2d}\times \mbb S$), we have $\widetilde{\mathcal{A}}f(x,k,y,l)=\mathcal{A}g(x,k)$.
	
	Throughout this section, we make the assumption that $\sup_{k\in\S}L_k<\infty$ as follows. As in the proof of Proposition 5.2.13 in \cite{1994Brownian}, for every fixed $k\in \S$, we can construct a sequence $\{\psi_{n}(r)\}_{n=1}^{\infty}$ of twice continuously differentiable functions satisfying $|\psi_{n}^\prime(r)|\leq1$ and $\lim_{n\to\infty}\psi_{n}(r)=|r|$ for $r\in \R$, and $0\leq\psi_{n}^{\prime\prime}(r)\leq\frac{2}{\sup_{k\in\S}L_knr^2}$ for $r\neq0$, where $L_k$ is as in Assumption \ref{asn:basic}. Furthermore, for every $r\in \R$, the sequence $\{\psi_{n}(r)\}_{n=1}^{\infty}$ is nondecreasing.
	
	\begin{lem}\label{lem:fellerpro}
		Under Assumptions \ref{asn:zysl}, \ref{asn:basic} and \ref{asn:feller}, for each $n\in  \mbb N$, let the function $\psi_{n}$ be defined as above and further define the function, for any $(x,k,y,l)\in \mbb R^{2d}\times \mbb S\times \mbb R^{2d}\times \mbb S $,
		$$f_{n}(x,k,y,l):=\psi_{n}(|x-y|)+\mathbf1_{\{k\neq l\}}.$$
		Then for all $(x,k,y,l)\in \mbb R^{2d}\times \mbb S\times \mbb R^{2d}\times \mbb S $ with $x\neq y$, we have
		\begin{equation}\label{eq:fellerest}
			\widetilde{\mathcal{A}}f_{n}(x,k,y,k)\leq\frac{1}{n}+C|x-y|,
		\end{equation}
		where $C=C(L_k,M_k,\Pi(U))$ is a positive constant.
	\end{lem}
	\begin{thm}\label{thm:feller}
		Under Assumptions \ref{asn:zysl}, \ref{asn:basic} and \ref{asn:feller}, the process $(X(t),\Lambda(t))$ generated by the operator $\mathcal{A}$ has the Feller property.
	\end{thm}
	\noindent\textbf{\textit{Proof.}} Denote by $\{P(t,(x,k),A):t\geq0,(x,k)\in \mbb R^{2d}\times \mbb S,A\in\mathcal B(\mbb R^{2d}\times \mbb S)\}$ the transition probability family of the process $(X(t),\Lambda(t))$. Since $\mbb S$ has a discrete topology, we only need to prove that for each $t\geq0$, $x$, $y\in \mbb  R^{2d} $ and $k\in \mbb S$, $P(t,(x,k),\cdot)$ converges weakly to $P(t,(y,k),\cdot)$ as $|x-y|\rightarrow0$. By virtue of Theorem 5.6 in \cite{2004From}, it suffices to prove that
	\begin{equation}\label{eq:wametric}
		W(P(t,(x,k),\cdot),P(t,(y,k),\cdot))\rightarrow0\quad {\rm as} \quad x\rightarrow y,
	\end{equation}
	where $W(\cdot,\cdot)$ denotes the Wasserstein metric between two probability measures.
	
	We can now utilize the aforementioned coupling to establish \eqref{eq:wametric}. Let $(\widetilde{X}(t),\widetilde\Lambda(t),\widetilde{Y}(t),\widetilde\Xi(t))$ denote the coupling process corresponding to the coupling operator $\widetilde{A}$ defined in \eqref{eq:couplingall}. Let us assume that $(\widetilde{X}(0),\widetilde\Lambda(0),\widetilde{Y}(0),\widetilde\Xi(0))=(x,k,y,k)\in \mbb R^{2d}\times \mbb S\times \mbb R^{2d}\times \mbb S $ with $x\neq y$. Define $\zeta:={\rm inf }\{t\geq0:\widetilde\Lambda(t)\neq\widetilde\Xi(t)\}$. Note that $\mbb P\{\zeta>0\}=1$. As in the proof of Theorem 2.3 in \cite{1989Coupling}, let us set 
	$$T_{R}:={\rm inf} \{t\geq0:|\widetilde{X}_{1}(t)|^2+|\widetilde{X}_{2}(t)|^2+|\widetilde{Y}_{1}(t)|^2+|\widetilde{Y}_{2}(t)|^2+\widetilde\Lambda(t)+\widetilde\Xi(t)>R\}.$$
	Using It\^o's formula to the process $f_{n}(\widetilde{X}(\cdot),\widetilde\Lambda(\cdot),\widetilde{Y}(\cdot),\widetilde\Xi(\cdot))$, we have 
	\begin{equation}\label{eq:fellerchang}
		\begin{aligned}
			&\mbb E\left[f_{n}(\widetilde{X}(t\wedge T_{R}\wedge\zeta),\widetilde\Lambda(t\wedge T_{R}\wedge\zeta),\widetilde{Y}(t\wedge T_{R}\wedge\zeta),\widetilde\Xi(t\wedge T_{R}\wedge\zeta))\right]\\
			=&f_{n}(x,k,y,k)+\mbb E\left[\int_{0}^{t\wedge T_{R}\wedge\zeta}\widetilde{\mathcal{A}}f_{n}(\widetilde{X}(s),\widetilde\Lambda(s),\widetilde{Y}(s),\widetilde\Xi(s))\text{d}s\right]\\
			\leq&\psi_{n}(|x-y|)+\frac{t}{n}+C\mbb E\left[\int_{0}^{t\wedge T_{R}\wedge\zeta}|\widetilde{X}(s)-\widetilde{Y}(s)|\text{d}s\right]
		\end{aligned}
	\end{equation}
	where the last step follows from the observation that $\widetilde\Lambda(s)=\widetilde\Xi(s)=k$ for all $s\in[0,t\wedge T_{R}\wedge\zeta)$ and the estimate in \eqref{eq:fellerest}. In addition, since $f_{n}(x,k,y,l)=\psi_{n}(|x-y|)+\mathbf1_{\{k\neq l\}}\geq\psi_{n}(|x-y|)$, we have from \eqref{eq:fellerchang} that
	\begin{equation*}
		\mbb E\left[\psi_{n}(|\widetilde{X}(t\wedge T_{R}\wedge\zeta)-\widetilde{Y}(t\wedge T_{R}\wedge\zeta)|)\right]\leq\psi_{n}(|x-y|)+\frac{t}{n}+C\mbb E\left[\int_{0}^{t\wedge T_{R}\wedge\zeta}|\widetilde{X}(s)-\widetilde{Y}(s)|\text{d}s\right]
	\end{equation*}
	Recall that $\psi_{n}(|x|)\uparrow|x|$ as $n\rightarrow\infty$. Therefore, passing to the limit as $n\rightarrow\infty$ on both sides of the above equation, it follows from the Monotone Convergence Theorem that
	\begin{equation*}
		\begin{aligned}
			\mbb E\left[|\widetilde{X}(t\wedge T_{R}\wedge\zeta)-\widetilde{Y}(t\wedge T_{R}\wedge\zeta)|\right]&\leq|x-y|+C\mbb E\left[\int_{0}^{t\wedge T_{R}\wedge\zeta}|\widetilde{X}(s)-\widetilde{Y}(s)|\text{d}s\right]\\
			&=|x-y|+C\mbb E\left[\int_{0}^{t}|\widetilde{X}(s\wedge T_{R}\wedge\zeta)-\widetilde{Y}(s\wedge T_{R}\wedge\zeta)|\text{d}s\right]
		\end{aligned}
	\end{equation*}
	Thus, by the Gronwall's inequality, we obtain
	$$\mbb E\left[|\widetilde{X}(t\wedge T_{R}\wedge\zeta)-\widetilde{Y}(t\wedge T_{R}\wedge\zeta)|\right]\leq|x-y|\e^{Ct}$$
	Finally, letting $R\to\infty$, we get that
	\begin{equation}\label{eq:fellerineq}
		\mbb E\left[|\widetilde{X}(t\wedge\zeta)-\widetilde{Y}(t\wedge\zeta)|\right]\leq|x-y|\e^{Ct}.
	\end{equation}
	
	On the other hand, note that $\zeta\leq t$ if and only if $\widetilde{\Lambda}(t\wedge\zeta)\neq\widetilde{\Xi}(t\wedge\zeta)$. Setting $f(x,k,y,l):=\mathbf1_{\{k\neq l\}}$ and applying It\^o's formula formula to the process $f(\widetilde{X}(t),\widetilde\Lambda(t),\widetilde{Y}(t),\widetilde\Xi(t))$ implies 
	\begin{equation}\label{eq:felleranother}
		\begin{aligned}
			\mbb P\{\zeta\leq t\}&=\mbb E[\mathbf1_{\{\widetilde{\Lambda}(t\wedge\zeta)\neq\widetilde{\Xi}(t\wedge\zeta)\}}]=\mbb E[f(\widetilde{X}(t\wedge\zeta),\widetilde\Lambda(t\wedge\zeta),\widetilde{Y}(t\wedge\zeta),\widetilde\Xi(t\wedge\zeta))]\\
			&=\mbb E\left[\int_{0}^{t\wedge\zeta}\widetilde{\Omega}_sf(\widetilde{X}(s),\widetilde\Lambda(s),\widetilde{Y}(s),\widetilde\Xi(s))\text{d}s\right]\leq M_k \mbb E\left[\int_{0}^{t\wedge\zeta}|\widetilde{X}(s)-\widetilde{Y}(s)|\text{d}s\right]\\
			&\leq \frac{M_k}{C}|x-y|\e^{Ct},
		\end{aligned}
	\end{equation}
	where the first inequality above follows from \eqref{eq:fellerswitch} and the last step follows from \eqref{eq:fellerineq}. 
	
	Finally, note that $|\widetilde{X}(t)-\widetilde{Y}(t)|$ is integrable. Hence, 
	\begin{equation}\label{eq:fellerfinally}
		\begin{aligned}
			\mbb E[|\widetilde{X}(t)-\widetilde{Y}(t)|]&=\mbb E[|\widetilde{X}(t)-\widetilde{Y}(t)|\mathbf1_{\{\zeta\leq t\}}]+\mbb E[|\widetilde{X}(t)-\widetilde{Y}(t)|\mathbf1_{\{\zeta>t\}}]\\
			&\leq \mbb E[|\widetilde{X}(t)-\widetilde{Y}(t)|\mathbf1_{\{\zeta\leq t\}}]+\mbb E[|\widetilde{X}(t\wedge\zeta)-\widetilde{Y}(t\wedge\zeta)|],
		\end{aligned}
	\end{equation}
	which, togther with \eqref{eq:fellerineq} and \eqref{eq:felleranother}, implies that $\mbb E[|\widetilde{X}(t)-\widetilde{Y}(t)|]$ converges to $0$ as $x\rightarrow y$. 
	
	On the other hand, observe that if $\widetilde{\Lambda}(t)\neq\widetilde{\Xi}(t)$, then $\zeta\leq t$. Thus thanks to \eqref{eq:felleranother}, we also have 
	\begin{equation}\label{eq:fellerfffff}
		\mathbb{E}\left[\mathbf{1}_{\{\widetilde{\Lambda}(t) \neq \widetilde{\Xi}(t)\}}\right] \leq \mathbb{P}\{\zeta \leq t\} \leq \frac{M_k}{C}|x-y| \e^{C t}.
	\end{equation}
	Recall that $\lambda((x, m),(y, n)):=|x-y|+\mathbf{1}_{\{m \neq n\}}$ is a metric on $\mathbb{R}^{2d} \times \mathbb{S}$. Hence, combining \eqref{eq:fellerfinally} and \eqref{eq:fellerfffff} yields 
	$$
	\mathbb{E}[\lambda((\widetilde{X}(t), \widetilde{\Lambda}(t)),(\widetilde{Y}(t), \widetilde{\Xi}(t)))]=0,\quad {\rm as}\quad x\to y.
	$$
	This implies \eqref{eq:wametric} and therefore completes the proof.\hfill $\square$\par
	
	In the rest of this section, we will provide a proof for Lemma \ref{lem:fellerpro}.
	
	\noindent\textbf{\textit{Proof of Lemma \ref{lem:fellerpro}.}} For any $x$, $y\in \mbb R^d$ and $k$, $l\in \mbb S$, set
	%\begin{center}
	%begin{tabular}
	\[
	\begin{aligned}
		&A(x,k,y,l)=a(x,k)+a(y,l)-2\sigma(x,k)\sigma(y,l)^\prime,\\
		&\widehat{B}(x,k,y,l)=\langle x_{2}-y_{2},b(x,k)-b(y,l)\rangle,
	\end{aligned}
	\]
	%\end{tabular}
	%\end{center}
	and
	$$\overline{A}(x,k,y,l)=\langle x_{2}-y_{2},A(x,k,y,l)(x_{2}-y_{2})\rangle/|x_{2}-y_{2}|^2.$$
	Then as in the proof of Theorem 3.1 in \cite{1989Coupling}, we can verify that
	\begin{equation*}
		\begin{aligned}
			\widetilde{\Omega}_{d}f_{n}(x,k,y,k)=&\frac{1}{2}\frac{\psi_{n}^\prime(|x-y|)}{|x-y|}[\text{tr}(A(x,k,y,k))-\overline{A}(x,k,y,k)+2\widehat{B}(x,k,y,k)+2\langle x_{1}-y_{1},x_{2}-y_{2}\rangle]\\
			&+\frac{1}{2}\psi_{n}^{\prime\prime}(|x-y|)\overline{A}(x,k,y,k)
		\end{aligned}
	\end{equation*}
	Observe that $\text{tr}(A(x,k,y,k))=\|\sigma(x,k)-\sigma(y,k)\|^2$ and
	\[
	\begin{aligned}
		&\overline{A}(x,k,y,k)=\frac{\langle x_{2}-y_{2},(\sigma(x,k)-\sigma(y,k))(\sigma(x,k)-\sigma(y,k))^\prime(x_{2}-y_{2})\rangle}{|x_{2}-y_{2}|^2}\leq\|\sigma(x,k)-\sigma(y,k)\|^2.
	\end{aligned}
	\]
	Thus, it follows from \eqref{eq:asnbasic2} of Assumption \ref{asn:basic} that
	\begin{equation}\label{eq:fellerdiffu}
		\begin{aligned}
			%\begin{split}
			\widetilde{\Omega}_{d}f_{n}(x,k,y,k)&\leq\frac{1}{2|x-y|}(2\|\sigma(x,k)-\sigma(y,k)\|^2+|x_{2}-y_{2}|^2+|b(x,k)-b(y,k)|^2+|x-y|^2)\\
			&\quad+\frac{1}{\sup_{i\in\S}L_in|x-y|^2}\|\sigma(x,k)-\sigma(y,k)\|^2\\
			&\leq(L_k+1)|x-y|+\frac{1}{n},
		\end{aligned}
	\end{equation}
	where the first inequality follows from the construction of the function $\psi_{n}$.
	
	Next, since $|\psi_{n}^\prime|\leq1$, we can use \eqref{eq:asnbasic2} of Assumption \ref{asn:basic} to compute
	\begin{equation}\label{eq:fellerju}
		\begin{aligned}
			%\begin{split}
			\widetilde{\Omega}_{j}f_{n}(x,k,y,k)=&\int_{U}[f_{n}(x_1,x_2+c(x,k,u),k,y_1,y_2+c(y,k,u),k)-f_{n}(x,k,y,k)\Pi(\text{d}u)\\
			=&\int_{U}(\psi_{n}(\sqrt{|x_{1}-y_{1}|^2+|(x_{2}+c(x,k,u))-(y_{2}+c(y,k,u)|^2)})-\psi_{n}(|x-y|))\Pi(\text{d}u)\\
			\leq&\int_{U}|\sqrt{|x_{1}-y_{1}|^2+|(x_{2}+c(x,k,u))-(y_{2}+c(y,k,u))|^2}-|x-y||\Pi(\text{d}u)\\
			\leq&\int_{U}|c(x,k,u)-c(y,k,u)|\Pi(\text{d}u)\\
			\leq& \sqrt{\Pi(U)L_k}|x-y|,
			%\end{split}
		\end{aligned}
	\end{equation}
	where Höder's inequality has been used in the last inequality above.
	
	Finally, we estimate $\widetilde{\Omega}_{s}f_{n}(x,k,y,k)$. By \eqref{eq:couoperQ}, we have
	\begin{equation}\label{eq:fellerswitch}
		\begin{aligned}
			\widetilde{\Omega}_{s}f_{n}(x,k,y,k)&=\sum_{i\in \mbb S}[q_{ki}(x)-q_{ki}(y)]^+(\mathbf1_{\{i\neq k \}}-\mathbf1_{\{k\neq k \}})+\sum_{i\in \mbb S}[q_{ki}(y)-q_{ki}(x)]^+(\mathbf1_{\{i\neq k \}}-\mathbf1_{\{k\neq k \}})\\
			& =\sum_{i\in\mathbb{S} \backslash\{k\} }|q_{ki}(x)-q_{ki}(y)|\leq M_k|x-y|,
		\end{aligned}
	\end{equation}
	where we utilize \eqref{eq:asnfeller2} to obtain the last inequality.
	
	Now plug \eqref{eq:fellerdiffu}, \eqref{eq:fellerju}, \eqref{eq:fellerswitch} into \eqref{eq:couplingall} yields \eqref{eq:fellerest}. This completes the proof.\hfill $\square$\par
	
	\section{Strong Feller Property}
	In order to prove the strong Feller property, we further make the following assumption:
	\begin{asn}\label{asn:sfyuanshi}
		For each $k\in \mbb S$, $X^{(k),0}(t)= (X^{(k),0}_1 (t),X^{(k),0}_2 (t))$ has the strong Feller property and has a positive transition probability density with respect to the Lebesgue measure.
	\end{asn}
	Then, we now introduce a family of jump diffusions under Assumption \ref{asn:sfyuanshi}. For each $k\in \mbb S$, let the jump diffusion $X^{(k)}(t)=(X_{1}^{(k)}(t),X_{2}^{(k)}(t))$ satisfy the following stochastic differential-integral equation:
	\begin{equation}\label{sfeq:guding}
		\left\{
		\begin{aligned}
			\text{d}X_{1}^{(k)}(t)&=X_{2}^{(k)}(t)\text{d}t,\\
			\text{d}X_{2}^{(k)}(t)&=b(X_{1}^{(k)}(t),X_{2}^{(k)}(t),k)\text{d}t+\sigma(X_{1}^{(k)}(t),X_{2}^{(k)}(t),k)\text{d}B(t)+\int_{U}c(X_{1}^{(k)}(t-),X_{2}^{(k)}(t-),k,u)N(\text{d}t,\text{d}u).
		\end{aligned}
		\right.
	\end{equation}
	\begin{lem}\label{lem:positivetg}
		Suppose that Assumption \ref{asn:sfyuanshi} holds. For each given $k\in \mbb S$, the jump-diffusion process $X^{(k)}(t)$ has the strong Feller property with a positive transition probability density with respect to the Lebesgue measure.
	\end{lem}
	
	To proceed, we first consider the strong Feller property for a special type of switching jump diffusion $(V(t),\Psi(t))$. Let the component $V(t)$ satisfy
	\begin{equation}\label{eq:steq:fuzhu}
		\left\{
		\begin{aligned}
			\text{d}V_{1}(t)&=V_{2}(t)\text{d}t,\\
			\text{d}V_{2}(t)&=b(V_{1}(t),V_{2}(t),\Psi(t))\text{d}t+\sigma(V_{1}(t),V_{2}(t),\Psi(t))\text{d}B(t)+\int_{U}c(V_{1}(t-),V_{2}(t-),\Psi(t-),u)N(\text{d}t,\text{d}u),
		\end{aligned}
		\right.
	\end{equation}
	and the component $\Psi(t)$ that is independent of the Brownian motion $B(\cdot)$ and Poisson random measure $N(\cdot,\cdot)$  be a time-homogeneous Markov chain with state space $\mbb S$ satisfying
	\begin{equation}
		\mbb P\{\Psi(t+\Delta)=l|\Psi(t)=k\}=\left\{
		\begin{aligned}
			&\widehat{q}_{kl}\Delta+o(\Delta),  &&\mbox{if $l\neq k$},\\
			&1+\widehat{q}_{kk}\Delta+o(\Delta),   &&\mbox{if $l=k$}
		\end{aligned}
		\right.
	\end{equation}
	provided $\Delta\downarrow0$, where $\widehat{Q}=(\widehat{q}_{kl})$ is a conservation $Q$-matrix such that 
	\(\widehat{q}_{kl}=\sup_{x\in \R^{2d}}q_{kl}(x) \) for $k\neq l$ and \( \widehat{q}_{kk}=-\sum_{l\neq k} \widehat{q}_{kl}\) for all $k\in \S$. To ensure that $\widehat{Q}$ is well defined, we need to substitute Assumption \ref{asn:zysl} with the following assumption.
	\begin{asn}\label{asn:kunakun}
		Assume that for all $k\in \mbb S$, we have 
		\begin{equation}
			\widehat{q}_k:=-\widehat{q}_{kk}=\sum_{l\in \mbb S\backslash\{k\}}\sup_{x\in \R^{2d}}q_{kl}(x)\leq H
		\end{equation}
		for some constant $H>0$.
	\end{asn}
	In the sequel, we sometimes emphasize the process $(V(t), \Psi(t))$ with initial condition $(V(0), \Psi(0))=(x, k)$ by $\left(V^{(x, k)}(t), \Psi^{(k)}(t)\right)$. We denote by $\Gamma(t,(x,k),\cdot)$ the transition probability of $(V(t),\Psi(t))$. For subsequent use, let us fix a probability measure $\mu(\cdot)$ that is equivalent to the product measure on $\mbb R^{2d}\times \mbb S$ of the Lebesgue measure on $\mbb R^{2d}$ and the counting measure on $\mbb S$. 
	\begin{lem}\label{strlem:feller}
		Suppose that Assumptions \ref{asn:sfyuanshi} and \ref{asn:kunakun} hold. Then $(V(t),\Psi(t))$ has the strong Feller property and the transition probability $\Gamma(t,(x,k),\cdot)$ of $(V(t),\Psi(t))$ has density $\gamma(t,(x,k),\cdot)$ with respect to $\mu(\cdot)$.
	\end{lem}
	To derive Lemma \ref{strlem:basic}, Assumption \ref{asn:basic} needs to be replaced with
	\begin{asn}\label{asn:basic2}
		Suppose that for all $x,y\in \mbb R^{2d}$ and  $k\in \mbb S$, we have
		\begin{equation}\label{eq:asnbasic11}
			|b(x, k)|^{2}+\|\sigma(x, k)\|^{2}+\int_{U}|c(x, k, u)|^{2} \Pi(\mathrm{d} u) \leq L\left(1+|x|^{2}\right),
		\end{equation}
		\begin{equation}\label{eq:asnbasic21} 
			|b(x,k)-b(y,k)|^2+\|\sigma(x,k)-\sigma(y,k)\|^2+\int_{U}|c(x,k,u)-c(y,k,u)|^2\Pi(\d u)\leq L|x-y|^2,			
		\end{equation}
		where the constant $L$ is positive.
	\end{asn}
	\begin{lem}\label{strlem:basic}
		Suppose that Assumptions \ref{asn:kunakun} and \ref{asn:basic2} hold. Then, for all $T>0$, $\delta>0$ and $k\in \mbb S$, we have
		\begin{equation}
			\mbb P\left\{\sup\limits_{0\leq t\leq T}|V^{(x,k)}(t)-V^{(y,k)}(t)|\geq\delta\right\}\rightarrow0
		\end{equation}
		as $x$, $y\in \mbb  R^{2d}$, $|x-y|\rightarrow0$.
	\end{lem}
	\begin{lem}\label{strlem:prob}
		Suppose that Assumptions \ref{asn:sfyuanshi} and \ref{asn:kunakun} hold. For any $t>0$ and bounded measurable function $f$ on $\mbb R^{2d}\times \mbb S$, we have 
		\begin{equation}
			f(V^{(x,k)}(t),\Psi^{(k)}(t))\rightarrow f(V^{(y,k)}(t),\Psi^{(k)}(t))\qquad \text{in probability}
		\end{equation}
		as $|x-y|\rightarrow0$.
	\end{lem}
	In order to transfer the strong Feller property from $(V(t),,\Psi(t))$ to $(X(t),\Lambda(t))$, we need to make a comparison between these two processes. Let $\{v_{m}\}$ be the sequence of stopping times defined by
	$$v_{0}=0,\qquad v_{m+1}=\inf\{s>v_{m}:\Psi(t)\neq\Psi(v_{m})\}\quad \text{for}\quad m\geq0.$$
	Define $n(t)=\max\{m : v_{m}\leq t\}$, which is the number of switches of $\Psi$ up to time $t$. Set $D:=D([0,\infty),\mbb R^{2d}\times \mbb S)$ and denote by $\mathcal{D}$ the usual $\sigma$-field of $D$. Likewise, for any $T>0$, set $D_{T}:=D([0,T],\mbb R^{2d}\times \mbb S)$ and denote by $\mathcal{D}_{T}$ the usual $\sigma$-field of $D_{T}$. Moreover, denote by $\mu_{1}(\cdot)$ the probability distribution induced by $(X(t),\Lambda(t))$ and $\mu_{2}(\cdot)$ the probability distribution induced by $(V(t),\Psi(t))$ in the path space $(D,\mathcal{D})$, respectively. Denote by $\mu_{1}^{T}(\cdot)$ the restriction of $\mu_{1}(\cdot)$ and $\mu_{2}^{T}(\cdot)$ the restriction of $\mu_{2}(\cdot)$ to $(D_{T},\mathcal{D}_{T})$, respectively. For any given $T>0$, from Lemma 4.2 of \cite{2009Asymptotic}, we know that $\mu_{1}^T(\cdot)$ is absolutely continuous with respect to $\mu_{2}^T(\cdot)$ and the corresponding Radon-Nikodym derivative has the following form.
	\begin{equation}\label{streq:rnde}
		\begin{aligned}
			M_{T}(V(\cdot),\Psi(\cdot)):&=\frac{\text{d}\mu_{1}^T}{\text{d}\mu_{2}^T}(V(\cdot),\Psi(\cdot))\\
			&=\prod_{i=0}^{n(T)-1}\frac{q_{\Psi(v_{i})\Psi(v_{i+1})}(V(v_{i+1}))}{\widehat{q}_{\Psi(v_{i})\Psi(v_{i+1})}}\exp\Big(-\sum\limits_{i=0}^{n(T)}\int_{v_{i}}^{v_{i+1}\wedge T}[q_{\Psi(v_{i})}(V(s))-\widehat{q}_{\Psi(v_{i})}]\text{d}s\Big).
		\end{aligned}
	\end{equation}
	\begin{lem}\label{strlem:conver}
		Suppose that Assumption \ref{asn:kunakun} and \ref{asn:basic2} hold, for all $T>0$, we have that
		\begin{equation}
			\mbb E[|M_{T}(V^{(x,k)}(\cdot),\Psi^{(k)}(\cdot))-M_{T}(V^{(y,k)}(\cdot),\Psi^{(k)}(\cdot))|]\rightarrow0
		\end{equation}
		as $|x-y|\rightarrow0$.
	\end{lem}
	\begin{lem}\label{strlem:inte}
		Suppose that Assumption \ref{asn:kunakun} holds, for all $T>0$ and $(x,k)\in \mbb R^{2d}\times \mbb S$, $M_{T}(V^{(x,k)}(\cdot),\Psi^{(k)}(\cdot))$ is integrable.
	\end{lem}
	We now present the main conclusion of this section.
	\begin{thm}\label{thm:strongfeller}
		Suppose that Assumptions \ref{asn:sfyuanshi}, \ref{asn:kunakun} and \ref{asn:basic2} hold. Then $(X(t),\Lambda(t))$ has the strong Feller property.
	\end{thm}
	\noindent\textbf{\textit{Proof.}} To prove the desired strong Feller property, it is enough to prove that for any $t>0$ and any bounded measurable function $f$ on $\mbb R^{2d}\times \mbb S$, $ \mbb E[f(X_{1}^{(x,k)}(t),X_{2}^{(x,k)}(t),\Lambda^{(x,k)}(t))]$ is bounded continuous in both $x\in \mbb R^{2d}$ and $k\in \mbb S$. Since $\mbb S$ has a discrete metric, it is sufficient to prove that
	\begin{equation}
		\left|\mbb E\left[f\left(X_{1}^{(x,k)}(t),X_{2}^{(x,k)}(t),\Lambda^{(x,k)}(t)\right)\right]-\mbb E\left[f\left(X_{1}^{(y,k)}(t),X_{2}^{(y,k)}(t),\Lambda^{(y,k)}(t)\right)\right]\right|\rightarrow0
	\end{equation}
	as $|x-y|\rightarrow0$. Indeed, by the Radon-Nikodym derivative, for all $(x,k)\in \mbb R^{2d}\times \mbb S$,
	\begin{equation}
		\begin{split}
			&\mbb E\left[f\left(X_{1}^{(x,k)}(t),X_{2}^{(x,k)}(t),\Lambda^{(x,k)}(t)\right)\right]\\
			=&\mbb E\left[f\left(V_{1}^{(x,k)}(t),V_{2}^{(x,k)}(t),\Psi^{(k)}(t)\right)\cdot M_{t}\left(V_{1}^{(x,k)}(t),V_{2}^{(x,k)}(t),\Psi^{(k)}(t)\right)\right]
		\end{split}
	\end{equation}
	Therefore, for any given $\varepsilon>0$, we have
	\begin{equation}\label{streq:longeq}
		\begin{aligned}
			&\quad\left|\mbb E\left[f\left(X_{1}^{(x,k)}(t),X_{2}^{(x,k)}(t),\Lambda^{(x,k)}(t)\right)\right]-\mbb E\left[f\left(X_{1}^{(y,k)}(t),X_{2}^{(x,k)}(t),\Lambda^{(y,k)}(t)\right)\right]\right|\\
			&\leq \mbb E\left|f\left(V_{1}^{(x,k)}(t),V_{2}^{(x,k)}(t),\Psi^{(k)}(t)\right)\cdot M_{t}\left(V_{1}^{(x,k)}(t),V_{2}^{(x,k)}(t),\Psi^{(k)}(t)\right)\right.\\
			&\left. \qquad-f\left(V_{1}^{(y,k)}(t),V_{2}^{(y,k)}(t),\Psi^{(k)}(t)\right)\cdot M_{t}\left(V_{1}^{(y,k)}(t),V_{2}^{(y,k)}(t),\Psi^{(k)}(t)\right)\right|\\
			&\leq\|f\| \mbb E\left|M_{t}\left(V_{1}^{(x,k)}(\cdot),V_{2}^{(x,k)}(\cdot),\Psi^{(k)}(\cdot)\right)-M_{t}\left(V_{1}^{(y,k)}(\cdot),V_{2}^{(y,k)}(\cdot),\Psi^{(k)}(\cdot)\right)\right|\\
			&+2\|f\| \mbb E\left[M_{t}\left(V_{1}^{(y,k)}(\cdot),V_{2}^{(y,k)}(\cdot),\Psi^{(k)}(\cdot)\right)\mathbf{1}_{\{|f\left(V_{1}^{(x,k)}(t),V_{2}^{(x,k)}(t),\Psi^{(k)}(t)\right)-f\left(V_{1}^{(y,k)}(t),V_{2}^{(y,k)}(t),\Psi^{(k)}(t)\right)|\geq\varepsilon\}}\right]\\
			&+\varepsilon \mbb E\left[M_{t}\left(V_{1}^{(y,k)}(\cdot),V_{2}^{(y,k)}(\cdot),\Psi^{(k)}(\cdot)\right)\right]\\
			&=\mathrm{(I)}+\mathrm{(II)}+\mathrm{(III)},
		\end{aligned}
	\end{equation}
	where $\|f\|:=\sup\{|f(x,k)|:(x,k)\in \mbb R^{2d}\times \mbb S\}$. From Lemma \ref{strlem:conver}, term $\mathrm{(I)}$ in \eqref{streq:longeq} tends to zero as $|x-y|\rightarrow0$. From Lemma \ref{strlem:prob} and \ref{strlem:inte}, we derive that term $\mathrm{(II)}$ in \eqref{streq:longeq}  also tends to zero as $|x-y|\rightarrow0$. Meanwhile, term $\mathrm{(III)}$ in \eqref{streq:longeq}  can be arbitrarily small since the multiplier $\varepsilon$ is arbitrary and $M_{t}(V_{1}^{(y,k)}(\cdot),V_{2}^{(y,k)}(\cdot),\Psi^{(k)}(\cdot))$ is integrable by Lemma \ref{strlem:inte}. The proof is complete.\hfill $\square$\par
	
	In the remainder of this section, we will prove Lemmas \ref{lem:positivetg}, \ref{strlem:feller}, and \ref{strlem:basic}-\ref{strlem:inte}.
	
	\noindent\textbf{\textit{Proof of Lemma \ref{lem:positivetg}.}} For a given $k\in \mbb S$, let us denote by $P^{(k)}(t,x,A)$ the transition probability for the process $X^{(k)}(t)$, and by $P^{(k),0}(t,x,A)$ the transition probability for the process $X^{(k)}(t)$. For any given $t>0$, $x\in \mbb R^{2d}$ and $A\in \mathcal{B}(\mbb R^{2d})$, we obtain the relation
	\begin{equation}\label{sfeq:jsf}
		\begin{aligned}
			P^{(k)}(t,x,A)=&\text{exp}\{-t\Pi(U)\}P^{(k),0}(t,x,A)\\
			&+\int_{0}^{t}\int\int_{U}\text{exp}\{-s_{1}\Pi(U)\}P^{(k),0}(s_{1},x,\text{d}y_{1}^{(1)}\times \text{d}y_{2}^{(1)})\Pi(\text{d}u_{1})\text{d}s_{1}\\
			&\qquad\qquad\quad\times P^{(k)}(t-s_{1},y_{1}^{(1)},y_{2}^{(1)}+c(y_{1}^{(1)},y_{2}^{(1)},k,u_{1}),A).
		\end{aligned}
	\end{equation}
	From this, we have
	\begin{equation*}
		\begin{aligned}
			&P^{(k)}(t-s_{1},y_{1}^{(1)},y_{2}^{(1)}+c(y_{1}^{(1)},y_{2}^{(1)},k,u_{1}),A)\\
			=&\text{exp}\{-(t-s_{1})\Pi(U)\}P^{(k),0}(t-s_{1},y_{1}^{(1)},y_{2}^{(1)}+c(y_{1}^{(1)},y_{2}^{(1)},k,u_{1}),A)\\
			&+\int_{0}^{t-s_{1}}\int\int_{U}\text{exp}\{-s_{2}\Pi(U)\}P^{(k),0}(s_{2},y_{1}^{(1)},y_{2}^{(1)}+c(y_{1}^{(1)},y_{2}^{(1)},k,u_{1}),\text{d}y_{1}^{(2)}\times \text{d}y_{2}^{(2)})\Pi(\text{d}u_{2})\text{d}s_{2}\\
			&\qquad\qquad\qquad\times P^{(k)}(t-s_{1}-s_{2},y_{1}^{(2)},y_{2}^{(2)}+c(y_{1}^{(2)},y_{2}^{(2)},k,u_{2}),A).
		\end{aligned}
	\end{equation*}
	Using \eqref{sfeq:jsf} again, we further have
	\begin{equation*}
		\begin{aligned}
			&P^{(k)}(t-s_{1}-s_{2},y_{1}^{(2)},y_{2}^{(2)}+c(y_{1}^{(2)},y_{2}^{(2)},k,u_{2}),A)\\
			=&\text{exp}\{-(t-s_{1}-s_{2})\Pi(U)\}P^{(k),0}(t-s_{1}-s_{2},y_{1}^{(2)},y_{2}^{(2)}+c(y_{1}^{(2)},y_{2}^{(2)},k,u_{2}),A)\\
			&+\int_{0}^{t-s_{1}-s_{2}}\int\int_{U}\text{exp}\{-s_{3}\Pi(U)\}P^{(k),0}(s_{3},y_{1}^{(2)},y_{2}^{(2)}+c(y_{1}^{(2)},y_{2}^{(2)},k,u_{2}),\text{d}y_{1}^{(3)}\times \text{d}y_{2}^{(3)})\Pi(\text{d}u_{3})\text{d}s_{3}\\
			&\qquad\qquad\quad\qquad\quad\times P^{(k)}(t-s_{1}-s_{2}-s_{3},y_{1}^{(3)},y_{2}^{(3)}+c(y_{1}^{(3)},y_{2}^{(3)},k,u_{3}),A).
		\end{aligned}
	\end{equation*}
	Using \eqref{sfeq:jsf} countably many times, we conclude that for any given $t>0$, $x\in \mbb R^{2d}$ and $A\in \mathcal{B}(\mbb R^{2d})$,
	\begin{equation}\label{sfeq:js}
		P^{(k)}(t,x,A)=\text{a\ series}.
	\end{equation}
	For this series, from the above equations, we derive that the first term (in which the process has no jump on $[0,t]$) is
	\begin{equation*}
		\text{exp}\{-t\Pi(U)\}P^{(k),0}(t,x,A).
	\end{equation*}
	The second term (in which the process has just one jump on $[0,t]$) is
	\begin{equation*}
		\begin{aligned}
			&\text{exp}\{-t\Pi(U)\}\int_{0}^{t}\int\int_{U}P^{(k),0}(s_{1},x,\text{d}y_{1}^{(1)}\times \text{d}y_{2}^{(1)})\Pi(\text{d}u_{1})\text{d}s_{1}\\
			&\qquad\qquad\quad\times P^{(k),0}(t-s_{1},y_{1}^{(1)},y_{2}^{(1)}+c(y_{1}^{(1)},y_{2}^{(1)},k,u_{1}),A).
		\end{aligned}
	\end{equation*}
	The third term (in which the process has just two jumps on $[0,t]$) is
	\begin{equation*}
		\begin{aligned}
			&\text{exp}\{-t\Pi(U)\}\int_{0}^{t}\int\int_{U}\int_{0}^{t-s_{1}}\int\int_{U}P^{(k),0}(s_{1},x,\text{d}y_{1}^{(1)}\times \text{d}y_{2}^{(1)})\Pi(\text{d}u_{1})\text{d}s_{1}\\
			&\qquad\qquad\times P^{(k),0}(s_{2},y_{1}^{(1)},y_{2}^{(1)}+c(y_{1}^{(1)},y_{2}^{(1)},k,u_{1}),\text{d}y_{1}^{(2)}\times \text{d}y_{2}^{(2)})\Pi(\text{d}u_{2})\text{d}s_{2}\\
			&\qquad\qquad\times P^{(k),0}(t-s_{1}-s_{2},y_{1}^{(2)},y_{2}^{(2)}+c(y_{1}^{(2)},y_{2}^{(2)},k,u_{2}),A),
		\end{aligned}
	\end{equation*}
	and moreover, the general term (in which the process has just $n$ jumps on $[0,t]$) is
	\begin{equation*}
		\begin{aligned}
			&\text{exp}\{-t\Pi(U)\}\int_{0}^{t}\int\int_{U}\int_{0}^{t-s_{1}}\int\int_{U}\cdots\int_{0}^{t-s_{1}-\cdots-s_{n-1}}\int\int_{U}P^{(k),0}(s_{1},x,\text{d}y_{1}^{(1)}\times \text{d}y_{2}^{(1)})\Pi(\text{d}u_{1})\text{d}s_{1}\\
			&\qquad\qquad\times P^{(k),0}(s_{2},y_{1}^{(1)},y_{2}^{(1)}+c(y_{1}^{(1)},y_{2}^{(1)},k,u_{1}),\text{d}y_{1}^{(2)}\times \text{d}y_{2}^{(2)})\Pi(\text{d}u_{2})\text{d}s_{2}\times\cdots\\
			&\qquad\qquad\times P^{(k),0}(s_{n},y_{1}^{(n-1)},y_{2}^{(n-1)}+c(y_{1}^{(n-1)},y_{2}^{(n-1)},k,u_{n-1}),\text{d}y_{1}^{(n)}\times \text{d}y_{2}^{(n)})\Pi(\text{d}u_{n})\text{d}s_{n}\\
			&\qquad\qquad\times P^{(k),0}(t-s_{1}-s_{2}-\cdots-s_{n},y_{1}^{(n)},y_{2}^{(n)}+c(y_{1}^{(n)},y_{2}^{(n)},k,u_{n}),A).
		\end{aligned}
	\end{equation*}
	In general, it is easy to see that the $n$-th term does not exceed
	$\frac{(t\Pi(U))^{n-1}}{(n-1)!}\text{exp}\{-t\Pi(U)\}.$
	Hence it follows that the series in \eqref{sfeq:js} converges uniformly with respect to $x$ over $\mbb R^{2d}$.
	
	It is easy to prove that for any given $t>0$ and $A\in \mathcal{B}(\mbb R^{2d})$, each term of the series in \eqref{sfeq:js} is lower semicontinuous with respect to $x$ by the strong Feller property of $X^{(k),0}(t)$. Therefore, it follows that for any given $t>0$ and $A\in \mathcal{B}(\mbb R^{2d})$, $P^{(k)}(t,x,A)$ is also lower semicontinuous with respect to $x$. As a result, $X^{(k)}(t)$ has the strong Feller property by Proposition 6.1.1 in \cite{1999Meyta}. Finally, from \eqref{sfeq:js}, $X^{(k)}(t)$ has a positive transition probability density with respect to the Lebesgue measure since $X^{(k),0}(t)$ does so under Assumption \ref{asn:sfyuanshi}. The proof is complete. \hfill $\square$\par
	
	\noindent\textbf{\textit{Proof of Lemma \ref{strlem:feller}. }} Denote by the $\nu_{1}$ the stopping time defined by $\nu_{1}=\inf\{s>0:\Psi(t)\neq\Psi(0)\}$. When $\Psi(0)=k$, $(\nu_{1},\Psi(\nu_{1}))$ on $[0,\infty)\times \mbb S_{k}$ with respect to the product of the Lebesgue measure and the counting measure has the probability density \text{exp}$(- \widehat{q}_k s)\widehat{q}_{kl}$, where $\mbb S_{k}:=\{l\in \mbb S:\widehat{q}_{kl}>0\}$. For any given $t>0$, $x\in \mbb R^{2d}$, $k$, $l\in \mbb S$ and $A\in \mcl B(\mbb R^{2d})$, we have the relation
	\begin{equation}\label{streq:count}
		\begin{aligned}
			\Gamma(t,(x,k),A\times\{l\})=&\delta_{kl}\text{exp}\{-\widehat{q}_k t\}P^{(k)}(t,x,A)\\
			&+\int_{0}^{t}\sum\limits_{l_{1}\in \mbb S_{k}}\int \widehat{q}_{kl_1}\text{exp}\{-\widehat{q}_{k} s_{1}\}P^{(k)}(s_{1},x,\text{d}y_{1}^{(1)}\times \text{d}y_{2}^{(1)})\\
			&\times\Gamma(t-s_{1},(y_{1}^{(1)},y_{2}^{(1)},l_{1}),A\times\{l\})\text{d}s_{1},
		\end{aligned}
	\end{equation}
	where $\delta_{kl}$ is the Kronecker symbol in $k$, $l$ which equals $1$ if $k=l$ and is $0$ if $k\neq l$. From this, we have
	\begin{equation*}
		\begin{aligned}
			\Gamma(t-s_{1},(y_{1}^{(1)},y_{2}^{(1)},l_{1}),A\times\{l\})=&\delta_{l_{1}l}\text{exp}\{-\widehat{q}_{l_1} (t-s_{1})\}P^{(l_{1})}(t-s_{1},y_{1}^{(1)},y_{2}^{(1)},A)\\
			&+\int_{0}^{t-s_{1}}\sum\limits_{l_{2}\in \mbb S_{l_{1}}}\int \widehat{q}_{l_1l_2}\text{exp}\{- \widehat{q}_{l_1}s_{2}\}P^{(l_{1})}(s_{2},y_{1}^{(1)},y_{2}^{(1)},\text{d}y_{1}^{(2)}\times \text{d}y_{2}^{(2)})\\
			&\times\Gamma(t-s_{1}-s_{2},(y_{1}^{(2)},y_{2}^{(2)},l_{2}),A\times\{l\})\text{d}s_{2}
		\end{aligned}
	\end{equation*}
	Using \eqref{streq:count} countably many times, as in the proof of Lemma \ref{lem:positivetg}, we conclude that for any given $t>0$, $x\in \mbb R^{2d}$ and $A\in \mcl{B}(\mbb R^{2d})$,
	\begin{equation}\label{streq:series}
		\Gamma(t,(x,k),A\times\{l\})=\text{a \ series}.
	\end{equation}
	For this series, we derive that the first term (in which $\Psi(t)$ has no jump on $[0,t]$) is
	\begin{equation*}
		\delta_{kl}\text{exp}\{- \widehat{q}_k t\}P^{(k)}(t,x,A).
	\end{equation*}
	The second term (in which $\Psi(t)$ has just one jump on $[0,t]$) is
	\begin{equation*}
		\delta_{l_{1}l}\int_{0}^{t}\sum\limits_{l_{1}\in \mbb S_{k}}\int \widehat{q}_{kl_1}\text{exp}\{-\widehat{q}_{k} s_{1}\}\text{exp}\{-\widehat{q}_{l_1} (t-s_{1})\}P^{(k)}(s_{1},x,\text{d}y_{1}^{(1)}\times \text{d}y_{2}^{(1)})P^{(l_{1})}(t-s_{1},y_{1}^{(1)},y_{2}^{(1)},A)\text{d}s_{1},
	\end{equation*}
	and the third term (in which $\Psi$ has just two jumps on $[0,t]$) is
	\begin{equation*}
		\begin{aligned}
			&\delta_{l_{2}l}\int_{0}^{t}\int_{0}^{t-s_{1}}\sum\limits_{l_{1}\in \mbb S_{k},l_{2}\in \mbb S_{l_{1}}}\int\int \widehat{q}_{kl_1}\text{exp}\{-\widehat{q}_{k} s_{1}\}\widehat{q}_{l_1l_2}\text{exp}\{-\widehat{q}_{l_1} s_2\}\text{exp}\{-\widehat{q}_{l_2} (t-s_{1}-s_2)\}\\
			&\quad\times P^{(k)}(s_{1},x,\text{d}y_{1}^{(1)}\times \text{d}y_{2}^{(1)})P^{(l_{1})}(s_{2},y_{1}^{(1)},y_{2}^{(1)},\text{d}y_{1}^{(2)}\times \text{d}y_{2}^{(2)})P^{(l_{2})}(t-s_{1}-s_{2},y_{1}^{(2)},y_{2}^{(2)},A)\text{d}s_{2}\text{d}s_{1}.
		\end{aligned}
	\end{equation*}
	Similar to the proof of Lemma \ref{lem:positivetg}, by Assumption \ref{asn:kunakun}, we can easily verify that the $n$-th term of the series in \eqref{streq:series} is bounded above by $\frac{(Ht)^{n-1}}{(n-1)!}$. Thus it is uniformly convergent with respect to $x\in \mbb R^{2d}$. Noting that $\mbb S$ is a infinitely countable set with a discrete metric, and using similar arguments as those in the proof of Lemma \ref{lem:positivetg}, we derive this lemma. \hfill $\square$\par

	\noindent\textbf{\textit{Proof of Lemma \ref{strlem:basic}-\ref{strlem:inte}.}}
	The proof of Lemma \ref{strlem:prob} is analogous to that of Lemma 5.6 in \cite{2017feller}. Lemma \ref{strlem:basic}, \ref{strlem:conver}, and \ref{strlem:inte} can be established using the methods presented in Lemma 4.1, 4.3, and 4.4 of \cite{2009Asymptotic}, respectively. We omit the specific details of these proofs.\hfill $\square$\par
	
	\section{Exponential Ergodicity}
	In this section, we will restrict our attention to the exponential ergodicity for $(X(t),\Lambda(t))$. As in \cite{1993Meyt3}, for any positive function $f(x, k) \geq 1$ defined on $\mathbb{R}^{2d} \times \mbb S$ and any signed measure $v(\cdot)$ defined on $\mathcal{B}\left(\mathbb{R}^{2d} \times \mbb S\right)$, we write
	\[
	\|v\|_{f}=\sup \{|v(g)|: \text{all measurable $g(x, k)$ satisfying $|g| \leq f$} \},
	\]
	where $v(g)$ denotes the integral of function $g$ with respect to measure $v$. Note that the total variation norm $\|v\|$ is just $\|v\|_{f}$ in the special case when $f \equiv 1$. For a function $1\leq f<\infty $ on $\mathbb{R}^{2d} \times \mbb S$, Markov process $(X(t), \Lambda(t))$ is said to be $f$-exponentially ergodic if there exist a probability measure $\pi(\cdot)$, a constant $\theta<1$ and a finite-valued function $\Theta(x, k)$ such that
	$$
	\|P(t,(x, k), \cdot)-\pi(\cdot)\|_{f} \leq \Theta(x, k) \theta^{t}
	$$
	for all $t \geq 0$ and all $(x, k) \in \mathbb{R}^{2d} \times \mbb S$. 
	
	\begin{asn}\label{asn:irreducible}
		Assume that the matrix $Q$ is irreducible on $\mathbb{R}^{2 d}$ in the following sense: for any distinct $k$, $l \in \mathbb{S}$, there exist $r \in \mathbb{N}$, $k_{0}, k_{1}, \ldots, k_{r} \in \mathbb{S}$ with $k_{i} \neq k_{i+1}$,  $k_{0}=k$ and $k_{r}=l$ such that the set $\left\{x \in \mathbb{R}^{2 d}: q_{k_{i} k_{i+1}}(x)>0\right\}$ has positive Lebesgue measure for $i=0$, $1, \ldots, r-1$.
	\end{asn}
	
	\begin{lem}\label{lem:skelton}
		Suppose that Assumptions  \ref{asn:zysl}, \ref{asn:sfyuanshi} and \ref{asn:irreducible} hold. Then $(X(t), \Lambda(t))$ is $\mu$-irreducible, where $\mu(\cdot)$ is a product measure on $\mathbb{R}^{2 d} \times \mathbb{S}$ of the Lebesgue measure on $\mathbb{R}^{2 d}$ and the counting measure on $\mathbb{S}$. Moreover, for any given $\delta>0$, all compact subsets of $\mathbb{R}^{2 d} \times \mathbb{S}$ are petite for the $\delta$-skeleton chain of $(X(t), \Lambda(t))$.
	\end{lem}
	
	Next, a nonnegative function $V(x, k)$ defined on $\mathbb{R}^{2d} \times \mbb S$ is called a norm-like function if $V(x, k) \rightarrow \infty$ as $|x|\vee k \rightarrow \infty$. Now we also need to introduce another Foster-Lyapunov drift condition as follows. For some $\alpha$, $\beta>0$ and a norm-like function $V(\cdot, k) \in C^{2}\left(\mathbb{R}^{2d}\right)$ with $k \in \mbb S$,
	\begin{equation}\label{eeeq:drift}
		\mathcal{A} V(x, k) \leq-\alpha V(x, k)+\beta, \quad(x, k) \in \mathbb{R}^{2d} \times \mbb S.
	\end{equation}
	\begin{thm}\label{thm:exerg}
		Suppose Assumptions \ref{asn:sfyuanshi}, \ref{asn:kunakun}, \ref{asn:basic2}, \ref{asn:irreducible} and \eqref{eeeq:drift} hold. Then the process $(X(t), \Lambda(t))$ is $f$-exponentially ergodic with $f(x, k)=V(x, k)+1$ and $\Theta(x, k)=B(V(x, k)+1)$, where $B$ is a finite constant.
	\end{thm}
	
	\noindent\textbf{\textit{Proof.}} For any given constant $\delta>0$, from Lemma \ref{lem:skelton}, all compact sets of the state space $\mathbb{R}^{2d} \times \mbb S$ are petite for the $\delta$-skeleton chain $(X(n\delta), \Lambda(n\delta))_{n \geq 0}$. Consequently, using \eqref{eeeq:drift} and applying Theorem $6.1$ in \cite{1993Meyt3} to the Markov process $(X(t), \Lambda(t))$, we obtain the desired result. The proof is complete.\hfill $\square$\par

	\noindent\textbf{\textit{Proof of Lemma \ref{lem:skelton}.}} Firstly, prove that the process $(X(t), \Lambda(t))$ is $\mu$-irreducible (refer to \cite{1992meytw1,1993Meytw2} for the detailed definition of $\mu$-irreducibility) in the sense that for any $t>0$, $(x, k) \in \mathbb{R}^{2d} \times \mathbb{S}$, $A \in \mathcal{B}\left(\mathbb{R}^{2d}\right)$ with positive Lebesgue measure and $l \in \mathbb{S}$, we have $P(t,(x, k), A \times\{l\})>0$. To do so, for each $k \in \mathbb{S}$, we kill  the process $X^{(k)}(t)$ of \eqref{sfeq:guding} with killing rate $q_{k}(\cdot)$. Denote by $\widetilde{P}^{(k)}(t, x, \cdot)$ the transition probability of the killed process. Then, by (5.13) in \cite{XI2021856} , we have
	\begin{equation*}
		\begin{split}
			&P(t,(x, k), A \times\{l\}) \\
			=&\delta_{k l} \widetilde{P}^{(k)}(t, x, A)+\sum_{m=1}^{\infty}~\idotsint\limits_{0<t_{1}<t_{2}<\cdots<t_{m}<t} \sum_{\substack{l_{0},l_{1}, l_{2}, \cdots, l_{m} \in \mathbb{S}\\l_{i} \neq l_ {i+1},l_{0}=k, l_{m}=l}} \int_{\mathbb{R}^{2d}} \cdots \int_{\mathbb{R}^{2d}} \widetilde{P}^{\left(l_{0}\right)}\left(t_{1}, x, \mathrm{d} z_{1}\right) q_{l_{0} l_{1}}\left(z_{1}\right) \\
			&\times \widetilde{P}^{\left(l_{1}\right)}\left(t_{2}-t_{1}, z_{1}, \mathrm{d} z_{2}\right) \cdots q_{l_{m-1} l_{m}}\left(z_{m}\right) \widetilde{P}^{\left(l_{m}\right)}\left(t-t_{m}, z_{m}, A\right) \mathrm{d} t_{1} \mathrm{d} t_{2} \cdots \mathrm{d} t_{m},
		\end{split}
	\end{equation*}
	where $\delta_{k l}$ is the Kronecker symbol, which equals 1 if $k=l$ and 0 if $k \neq l$. It is easy to see that for any $t>0$, $(x, k) \in \mathbb{R}^{2 d} \times \mathbb{S}$ and $A \subset \mathbb{R}^{2 d}$ with $A$ having positive Lebesgue measure,
	$$
	\widetilde{P}^{(k)}(t, x, A)=\mathbb{\mbb E}_{k}^{(x)}\left[\mathbf{1}_{A}\left(X^{(k)}(t)\right) \exp \left\{\int_{0}^{t} q_{k k}\left(X^{(k)}(u)\right) \mathrm{d} u\right\}\right] \geq \e^{-H t} P^{(k)}(t, z, A)>0
	$$
	where Assumption \ref{asn:zysl} and Lemma \ref{lem:positivetg} have been used. This, together with Assumption \ref{asn:irreducible}, implies that $P(t,(x, k), A \times\{l\})>0$, which is the desired assertion.
	
	Next, we show that for any given $\delta>0$, all compact subsets of $\mathbb{R}^{2 d} \times \mathbb{S}$ are petite for the $\delta$-skeleton chain  $(X(n\delta), \Lambda(n\delta))_{n \geq 0}$. By the $\mu$-irreducibility of the process $(X(t),\Lambda(t))$, its $\delta$-skeleton chain $(X(n \delta), \Lambda(n \delta))_{n \geq 0}$ are $\mu$-irreducible. Note
	that supp $\mu(\cdot)$ is equal to $\mathbb{R}^{2 d} \times \mathbb{S}$ and hence has non-empty interior. On the other hand, Theorem \ref{thm:strongfeller} says that $(X(t), \Lambda(t))$ is strong Feller and hence Feller. Combining these facts with Theorem 3.4 in \cite{1992meytw1}, we obtain that all compact subsets of $\mathbb{R}^{2 d} \times \mathbb{S}$ are petite for the $\delta$-skeleton chain of $(X(t), \Lambda(t))$. This completes the proof.\hfill $\square$\par
	
	\section{Example}	
	Let $\S=\{0, 1, 2,\ldots\}$ and $M$ be an arbitrary positive constant; let $\{\alpha(k) \}_{k\in\S}$ be a sequence satisfying $\alpha(k)<-2M^2$ and $\lim_{k \rightarrow \infty}\alpha(k)=-\infty$. For each $k\in \mbb S$, $\sigma(x_1, x_2, k)$ satisfies the Lipschitz condition such that $0<\sigma(x_1,x_2,k)\leq M$, and $\gamma(k)$ is any given real numbers. Consider the following stochastic process whose second component $X_{2}(t)$ is a coupled one-dimensional Ornstein-Uhlenbeck type process.  The diffusion process $X(t)=(X_{1}(t),X_{2}(t))$ satisfies the following stochastic differential equation in $\mbb R^{2}$,
	\begin{equation}
		\left\{
		\begin{aligned}
			\d X_{1}(t)&=X_{2}(t)\d t\\
			\d X_{2}(t)&=\alpha(\Lambda(t))X_{2}(t)\d t+\sigma(X_{1}(t),X_{2}(t),\Lambda(t))\d B(t)-\int_{\Gamma}X^{\frac{1}{3}}_2(t)\gamma(\Lambda(t-)) u N(\mathrm{d} t, \mathrm{d} u) ,
		\end{aligned}
		\right.
	\end{equation}
	where $B(t)$ is a standard one-dimensional Brownian motion and
	$N(\mathrm{d} t, \mathrm{d} u)$ is a stationary Poisson point process and independent of $B(t)$ such that $\widetilde{N}(\mathrm{d} t, \mathrm{d} u)=N(\mathrm{d} t, \mathrm{d} u)-\Pi(\mathrm{d} u) \mathrm{d} t$ is the compensated Poisson random measure on $[0, \infty) \times \mathbb{R}$, where $\Pi(\cdot)$ is a deterministic finite characteristic measure concentrated on the measurable space $(\Gamma,\mathcal{B}(\Gamma))$(here $\Gamma$ is a compact set not including the origin  0 in $\R$).
	
	Suppose the switching component $\Lambda(t)$ is generated by the $Q$-matrix, where ${Q}=({q}_{kl}(x))$ is a conservation $Q$-matrix such that for all $l\neq k$ satisfying $\alpha(l)<\alpha(k)$,
	\({q}_{kl}(x)=\exp(-2(\alpha(k)-\alpha(l))|x_1|) \) and $({q}_{kl}(x))$ sastifies Assumption \ref{asn:kunakun}.
	%		\[
	%		Q(x)=\begin{pmatrix}
		%			-q_{01}(x)&q_{01}(x)&0&0&0&0&\cdots\\
		%			q_{10}(x)&-(q_{10}(x)+q_{12}(x))&q_{12}(x)&0&0&0&\cdots\\
		%			0&q_{21}(x)&-(q_{21}(x)+q_{23}(x))&q_{23}(x)&0&0&\cdots\\
		%			\vdots&\vdots&\vdots&\vdots&\vdots&\vdots&\ddots
		%		\end{pmatrix},
	%		\]
	Apparently, Assumptions \ref{asn:sfyuanshi}, \ref{asn:kunakun} and \ref{asn:irreducible} hold.

	We need only show that  Condition \eqref{eeeq:drift} is satisfied.  In what follows, we verify that Condition \eqref{eeeq:drift} is also satisfied. To do so, as in the proof of \cite[Theorem 3.1]{2001Large}, we set a function ${V(x, k)}$ on $\mathbb{R}^{2} \times\S$ as
	\[
	V(x, k)=\exp \left(x_{2}^{2}+G(x_{1}) x_{2}+U\left(x_{1}, k\right)\right) .
	\]
	Here, the function ${G\left(x_{1}\right)}$ is infinitely differentiable in ${x_{1}}$ such that
	\begin{equation}\label{eq:sinian}
		G\left(x_{1}\right)=\frac{x_{1}}{\left|x_{1}\right|}  \text { for }\left|x_{1}\right|>1  \text { and } \left|G\left(x_{1}\right)\right| \leq 1 \text { for } x_{1} \in \mathbb{R} ;
	\end{equation}
	and the function ${U\left(x_{1}, k\right)}$ is twice differentiable in ${x_{1}}$ such that
	\begin{equation}\label{eq:sinianxp}
		U\left(x_{1}, k\right)=-\alpha(k)\left|x_{1}\right|,  \text { for }\left|x_{1}\right|>1  \text { and }  k \in\S .
	\end{equation}
	Clearly, ${V(x, k)}$ is a norm-like function. Moreover, for the operator ${\mathcal{A}}$ defined in \eqref{eq:alloperator},
	\[
	\begin{aligned}
		\mathcal{A} V(x, k)=& x_{2} \frac{\partial V(x, k)}{\partial x_{1}} +\alpha(k) x_{2} \frac{\partial V(x, k)}{\partial x_{2}} +\frac{1}{2} \sigma^{2}\left(x, k\right) \frac{\partial^{2}V(x, k)}{\partial x_{2}^{2}}  \\
		&+\int_{\Gamma}\left(V(x_1,x_2-x_2^{\frac{1}{3}}\gamma(k)u,k)-V(x_1,x_2,k) \right)\Pi(\d u)+\sum_{l \in\S \backslash\{k\}} q_{k l}(x)(V(x, l)-V(x, k)).\\
	\end{aligned}
	\]
	At the same time, by some elementary calculations, we also have
	\[
	\begin{aligned}
		\frac{\partial V(x, k)}{\partial x_{1}}&=V(x, k)\left[G^{\prime}\left(x_{1}\right) x_{2}+U^{\prime}\left(x_{1}, k\right)\right],\\
		\frac{\partial V(x, k)}{\partial x_{2}}&=V(x, k)\left[2 x_{2}+G\left(x_{1}\right)\right],\\
		\frac{\partial^{2} V(x, k)}{\partial x_{2}^{2}}&=V(x, k)\left[2 x_{2}+G\left(x_{1}\right)\right]^{2}+2 V(x, k),\\
		V(x_1,x_2-x_2^{\frac{1}{3}}\gamma(k)u,k)-V(x_1,x_2,k)&=V(x,k)\left[ \exp\left(x_2^{\frac{2}{3}}\gamma^2(k)u^2-(2x_2+G(x_1))x_2^{\frac{1}{3}}\gamma(k)u\right)-1\right]\\
		\sum_{l \in\S \backslash\{k\}} q_{k l}(x)(V(x, l)-V(x, k))&=V(x,k)\sum_{l \in\S \backslash\{k\}} q_{k l}(x)\Big(\exp\big((\alpha(k)-\alpha(l))|x_1|\big)-1\Big).
	\end{aligned}
	\]
	Inserting these and undergoing some tedious calculations, we then arrive at
	\begin{equation}\label{eq:wuliaoa}
		\mathcal{A} V(x, k)=-V(x, k) W_{k}(x),
	\end{equation}
	where
	\begin{equation}\label{eq:jieshu}
		\begin{aligned}
			W_{k}(x)=&-\frac{1}{2} \sigma^{2}\left(x, k\right)\left(2 x_{2}+G\left(x_{1}\right)\right)^{2}-\sigma^{2}\left(x, k\right) \\
			&-G^{\prime}\left(x_{1}\right) x_{2}^{2}-x_{2} U^{\prime}\left(x_{1}, k\right)-\alpha(k)x_{2} G\left(x_{1}\right)-2\alpha(k)x_2^2 \\
			&-\int_{\Gamma}\left[\exp\left(x_2^{\frac{2}{3}}\gamma^2(1)u^2-(2x_2+G(x_1))x_2^{\frac{1}{3}}\gamma(1)u \right)-1 \right]\Pi(\d u)\\
			&-\sum_{l \in\S \backslash\{k\}} q_{k l}(x)\Big(\exp\big((\alpha(k)-\alpha(l))|x_1|\big)-1\Big).
		\end{aligned}
	\end{equation}
	Let ${W_{k i}(x)}$ with ${i=1,2,3,4}$ be the ${i}$ th line on the right-hand side of \eqref{eq:jieshu}. Now let us estimate these functions one by one. First, noting ${0<\sigma\left(x, k\right) \leq M}$ and ${\left|G\left(x_{1}\right)\right| \leq 1}$, ${0 \geq W_{k1}(x) \geq-4 M^{2} x_{2}^{2}-2 M^{2}}$. Next, by virtue of \eqref{eq:sinian} and \eqref{eq:sinianxp}, for ${\left|x_{1}\right|>1}$,
	\[
	G^{\prime}\left(x_{1}\right)=0  \text { and }  U^{\prime}\left(x_{1}, k\right)=-\alpha(k) G\left(x_{1}\right).
	\]
	What's more, using \eqref{eq:sinian} , for ${\left|x_{1}\right|>1}$, $W_{k2}(x)=-2\alpha(k)x_2^2$. Moreover, for ${\left|x_{1}\right|>1}$,
	\[
	W_{k3}(x)=\Pi(\Gamma)-\int_{\Gamma}\exp\left(x_2^{\frac{2}{3}} \gamma^2(k)u^2-2x_2^{\frac{4}{3}}\gamma(k)u-\frac{x_1}{|x_1|}x_2^{\frac{1}{3}}\gamma(k)u \right)\Pi(\d u).
	\]
	Finally, using \eqref{eq:sinianxp} and recalling values of ${q_{kl}(x)}$ defined above, we get that for ${\left|x_{1}\right|>1}$,
	\[
	\begin{aligned}
		W_{k4}(x)=
		&\sum_{\substack{l \in\S \backslash\{k\}\\ \alpha(l)\geq\alpha(k)}} q_{k l}(x)\Big(1-\exp\big((\alpha(k)-\alpha(l))|x_1|\big)\Big)\\
		&+ \sum_{\substack{l \in\S \backslash\{k\}\\\alpha(l)<\alpha(k)}} \Big(\exp(-2(\alpha(k)-\alpha(l))|x_1|)-\exp(-(\alpha(k)-\alpha(l))|x_1|)\Big).
	\end{aligned}
	\]
	From the above estimations, it is easy to see that the leading term of the right-hand side of \eqref{eq:jieshu} is $-2(\alpha(k)+2M^2)x_2^2$. Hence, we can choose a compact subset ${C \subset \mathbb{R}^{2}}$ such that ${W_{k}(x) \geq 1}$ for all ${x \in C^{c}}$. Combining this with \eqref{eq:wuliaoa}, we obtain that
	\begin{equation}\label{eq:what}
		\begin{aligned}
			\mathcal{A} V(x, k) & \leq-V(x, k) \mathbf{1}_{C^{c}}(x)-V(x, k) W_{k}(x) \mathbf{1}_{C}(x) \\
			&=-V(x, k)+V(x, k)\left(1-W_{k}(x)\right) \mathbf{1}_{C}(x) \\
			& \leq-V(x, k)+\beta,
		\end{aligned}
	\end{equation}
	for some positive constant ${\beta>0}$.
	Therefore, by virtue of Theorem \ref{thm:exerg}, the strong Markov process ${(X(t), \Lambda(t))}$ is ${\Psi}$-exponentially ergodic with ${\Psi(x, k)=V(x, k)+1}$ and ${\Theta(x, k)=B(V(x, k)+1)}$, where ${B}$ is a finite constant.\qed
	\bibliographystyle{plain}
	\bibliography{math1}
\end{document}